\documentclass[reqno]{amsart}
\usepackage[english]{babel}
\usepackage{amscd,amssymb,amsmath,amsfonts,latexsym,amsthm}
\usepackage{inputenc}
\usepackage{graphicx}
\usepackage{hyperref,amsthm,amsmath,amssymb,color,multirow,graphicx}
\usepackage[fontsize=11pt]{scrextend}

\newtheorem{thm}{Theorem}[section]

\newtheorem{lem}[thm]{Lemma}

\theoremstyle{definition}
\newtheorem{defn}[thm]{Definition}

\numberwithin{equation}{section}

\newcommand{\be}{\begin{equation}}

\newcommand{\ee}{\end{equation}}

\begin{document}

\sloppy

\begin{flushright}
	\begin{tabular}{ll}
		\textsf{} &  \\
		\textsf{\ \pageref{firstpage}-\pageref{lastpage}}\\
		{}\\
\end{tabular}
\end{flushright}

\begin{center}

\textbf{\large Source recovery and finite-time blow-up for a generalized time-fractional Kuramoto-Sivashinsky equation}

 \vspace{5mm}

{\small
\begin{tabular}{p{9cm}}
    {\bf Sobirov Z.A.} \\
   National University of Uzbekistan, Tashkent, Uzbekistan,\\
   V.I. Romanovskiy Institute of Mathematics,\\
   Uzbekistan Academy of Sciences.\\
   sobirovzar@gmail.com\\
   \\[5mm]
    {\bf Norkulova R.B.}\\
   V.I. Romanovskiy Institute of Mathematics, \\
   Uzbekistan Academy of Sciences.\\
   norkulovarushana1988@gmail.com\\
    
\end{tabular}
}
\end{center}
\vspace{5mm}

\textbf{Abstract.} We study the recovery of an unknown time-dependent factor $f(t)$ in the source term $f(t)g(x,t)$ of a generalized Kuramoto--Sivashinsky equation with a Caputo time-fractional derivative of order $\beta\in(0,1)$. The problem is posed on $\mathbb{R}$ and is supplemented with an integral overdetermination condition. For the linear equation, the inverse problem is reduced to a second-kind operator equation, and the quasinilpotency of the corresponding operator is used to establish existence, uniqueness, and stability. For the nonlinear equation, a successive approximation scheme combined with a priori estimates yields local existence and uniqueness in the Schwartz space. Stability estimates with respect to perturbations of the initial data and the integral observation are also obtained. Finally, an explicit solution exhibiting finite-time blow-up is constructed, showing that the nondegeneracy condition on the integral observation alone does not guarantee global solvability.\\



\makeatletter
\renewcommand{\@evenhead}{\vbox{\thepage \hfil {\it  Sobirov Z.A., Norkulova R.B.}  \hrule }}
\renewcommand{\@oddhead}{\vbox{\hfill
{\it Sobirov Z.A., Norkulova R.B. }\hfill
\thepage \hrule}} \makeatother

\label{firstpage}

\section{INTRODUCTION AND STATEMENT OF THE PROBLEM}

The Kuramoto--Sivashinsky equation is a well-known nonlinear model that arises in fluid dynamics, plasma physics, and reaction--diffusion systems. It describes various nonlinear phenomena, including instabilities, turbulence, and pattern formation. The equation was independently derived by Kuramoto and Tsuzuki in the study of phase turbulence in reaction--diffusion systems and by Sivashinsky in the study of flame-front propagation \cite{Sivashinsky,Kuramoto,SIVASHINSKY}.

The classical Kuramoto--Sivashinsky equation takes the form
\[
u_t+u_{xxxx}+u_{xx}+uu_x=0.
\]
Here, the fourth-order dissipative term $u_{xxxx}$ suppresses energy at small scales, the second-order term $u_{xx}$ drives instability at large scales, and the nonlinear convective term $uu_x$ redistributes energy across scales.

In a direct problem, all coefficients, source terms, and initial or boundary conditions are assumed to be known, and the objective is to determine the solution of the governing equation. In many applications, however, some coefficients or source terms cannot be measured directly. Determining these unknown quantities from additional observations leads naturally to inverse problems. Such problems arise in engineering, geophysics, and other areas in which the parameters or sources governing a system are not directly observable.

Numerous authors have investigated direct problems for the Kuramoto--Sivashinsky equation and its fractional generalizations. In \cite{QASEM}, analytical solutions of a time-fractional Kuramoto--Sivashinsky equation were studied using the residual power series transform method based on the Caputo derivative, together with convergence analysis and numerical results. The authors \cite{Ashurov} proved the global solvability of the Cauchy problem for a generalized time-fractional Kuramoto--Sivashinsky equation in the Schwartz space, using the Fourier transform and the method of successive approximations. In this paper \cite{Irgashov}, the Cauchy problem in a half-plane was studied for a fourth-order inhomogeneous equation with a Caputo fractional derivative. Uniqueness was proved using the Laplace transform, while the Green function, an explicit solution, and a fundamental solution were constructed in terms of Wright functions.

The notions of inverse and ill-posed problems, together with numerous examples, are discussed in S.I.Kabanikhin \cite{Kabanikhin1}. A systematic theory to such problems for partial differential equations can be found in the monographs by S.I.Kabanikhin \cite{Kabanikhin 2011} and Prilepko et al. \cite{Prilepko}.

Inverse problems have been extensively investigated for various classes of differential equations. In particular, considerable attention has been devoted to source identification for diffusion and subdiffusion equations, fourth-order parabolic equations, and Kuramoto--Sivashinsky equations. Some representative works are briefly reviewed below.

For a time-fractional diffusion equation, the authors of \cite{Sakamoto} considered an inverse problem of recovering a spatially dependent source function from final overdetermination data and showed that the problem is well posed except for a discrete set of values of the diffusion coefficient. The researchers \cite{Ashurov1,Ashurov2} had studied the inverse problems of determining the  time-dependent or spatially dependent source function for the subdiffusion equation under different observation conditions. A space-time-dependent source identification problem for a subdiffusion equation with a Caputo derivative was investigated in \cite{Ashurov3}, where the existence and uniqueness of a weak solution were established using the Fourier method and the method of successive approximations. In these papers \cite{Sobirov1,Sobirov2,Sobirov3}, the inverse source problems for subdiffusion and space-time fractional equations on metric star graphs were reduced to operator equations, and the existence and continuity of the corresponding resolvent operators were established.
Inverse source problems for parabolic equations with integral overdetermination go back to the works of Prilepko \cite{Prilepko}, \cite{Prilepko1985}.

The unique solvability of an inverse problem for determining the right-hand side of a parabolic equation with a leading coefficient depending on time and space was investigated  \cite{Kamynin} under an integral overdetermination condition. Sufficient conditions for local solvability and Fredholm solvability were obtained. In \cite{Cao}, a space-dependent source term in a fourth-order parabolic equation was recovered from final integral measurements. The existence and uniqueness of the solution were established, and numerical methods for treating noisy data were proposed.

Inverse problems for the classical Kuramoto--Sivashinsky equation have been considered by several authors. In \cite{Baudouin}, the anti-diffusion coefficient in a one-dimensional nonlinear Kuramoto--Sivashinsky equation was recovered from solution measurements, and Lipschitz stability was established using the Bukhgeim--Klibanov method and a Carleman estimate. Bin Hu et al. \cite{BIN.} considered an ill-posed source inversion problem for a linear Kuramoto--Sivashinsky equation was regularized using the quasi-reversibility method. Convergence estimates were obtained for both a priori and a posteriori choices of the regularization parameter, and numerical examples were presented. K.Cao \cite{KAI.CAO} investigated  an unknown time-dependent source term in the Kuramoto--Sivashinsky equation was identified from integral observation data. Existence and uniqueness results were established, and a time-discrete reconstruction method was proposed with convergence estimates for exact and noisy data.

However, the above results do not directly cover the inverse source problem for a nonlinear generalized time-fractional
Kuramoto--Sivashinsky equation with a Caputo derivative. The memory effect of the fractional derivative, together
with the nonlinear convective term and the unknown time-dependent source, requires additional estimates for both the solution and the source function.

In the present paper, we study source recovery for a generalized time-fractional Kuramoto--Sivashinsky equation with a Caputo derivative. Specifically, we consider the equation
\begin{equation}\label{chiziqli tenglama}
\partial_{0t}^{\beta}u+a^2u_{xxxx}+bu_{xxx}+cu_{xx}+du_x+qu
=f(t)g(x,t)+\gamma u^pu_x,
\quad x\in\mathbb R,\quad 0<t\leq T,
\end{equation}
with the initial condition
\begin{equation}\label{3-tenglama}
u(x,0)=\varphi(x),\quad x\in\mathbb R,
\end{equation}
and the overdetermination condition
\begin{equation}\label{O.C}
\int\limits_{-\infty}^{+\infty}\eta(x)u(x,t)\,dx
=\psi(t),\quad 0\leq t\leq T.
\end{equation}
Here, $0<\beta<1$, $a>0$, $b,c,d,q,\gamma\in\mathbb R$, and $p$ is a positive integer. The operator $\partial_{0t}^{\beta}$ denotes the Caputo fractional derivative. The functions $g(x,t)$, $\eta(x)$, $\varphi(x)\in S(\mathbb R)$, and $\psi(t)$ are given, while the solution $u(x,t)$ and the time-dependent source factor $f(t)$ are unknown.

The additional condition \eqref{O.C} represents a weighted spatial integral of the state of the system and provides the observation used to recover the unknown source factor. The problem is studied in the Schwartz space $S(\mathbb R)$ of infinitely differentiable, rapidly decreasing functions.

For the linear equation, the inverse problem is reduced to a second-kind operator equation, and the quasinilpotency of the corresponding operator is used to prove existence, uniqueness, and stability. For the nonlinear problem, the method of successive approximations combined with appropriate a priori estimates yields local existence and uniqueness, with the state in the Schwartz space and the source factor continuous in time. Stability estimates with respect to perturbations of the initial data and the integral observation are also established. Finally, an explicit solution exhibiting finite-time blow-up is constructed. In this example, both the state and the source factor become unbounded, even though $\int_{\mathbb R}\eta(x)g(x,t)\,dx$ remains a strictly positive constant. Thus, the nondegeneracy condition alone does not guarantee global solvability.

The structure of this paper is as follows. Section 2 contains the necessary definitions and auxiliary results. In Section 3, we study the linear inverse problem and prove the existence, uniqueness, and stability of its solution. Section 4 is devoted to the local solvability of the nonlinear inverse problem. Stability estimates are established in Section 5. In Section 6, we construct an explicit solution exhibiting finite-time blow-up. The main conclusions are summarized in Section 7.

\section{AUXILIARY PROPOSITIONS}

In this section, we collect  the main results, facts, and notation that will be used in the remainder of the paper.
\subsection{The Schwartz space and the Fourier transform}
\begin{defn}\label{1def3}(\cite{John},\cite{Elias},p.135) The Schwartz space $S({{\mathbb R}})$ is the topological vector space of functions $g:{{\mathbb R}}\to \mathbb R $ such that $g\in {{C}^{\infty }}({{\mathbb R}})$ and
\begin{equation}\label{S(R)}\sup_{x\in \mathbb R}{|x|^{k}}\left|\frac{d^n g(x)}{dx^n}\right|<\infty , \end{equation}	
for every  $k,n\in \mathbb N\cup \{0\}$.\end{defn}

The Schwartz space is a natural one to use for the Fourier transform. Differentiations and multiplication exchange roles under the Fourier transform due to the properties of smoothness and rapid decrease. As a result, the Fourier transform is an automorphism of the Schwartz space.

\begin{defn}(\cite{John},\cite{Elias})
The Fourier transform of a function $f\in L^1({\mathbb R})$ is a function $\hat{f}:{\mathbb R}\to \mathbb C$ defined by  
\begin{equation}\hat{f}(\lambda )\equiv \mathbf F[f](\lambda):=\frac{1}{\sqrt{2\pi}}\int\limits_{-\infty }^{+\infty }{f(x){{e}^{-i\lambda x}}dx}. \end{equation}		
The inverse Fourier transform defined by 
\begin{equation}f(x)\equiv \mathbf F^{-1}[\hat f](x):=\frac{1}{\sqrt{2\pi} }\int\limits_{-\infty }^{+\infty }{\hat{f}(\lambda ){{e}^{i\lambda x}}d\lambda, }\end{equation} provided the right-hand side exists.

We  use $x$ to denote the spatial variable and $\lambda $  the variable in the Fourier transform.

\begin{thm}(\cite{Elias},p.137)
 \label{thm1}
The Fourier transform and its inverse map the space $S$ on to itself in a one-to-one, linear, and continuous manner. \end{thm}

\end{defn}

\subsection{Fractional integral and derivatives}

The Riemann-Liouville fractional integral of order $\beta>0$, of a function $g$ is determined by (see \cite{Kilbas2006}, p. 69)
$$I_{0t}^\beta g(t)=\frac{1}{\Gamma(\beta)}\int\limits_{0}^{t}\frac{g(\xi)}{(t-\xi)^{1-\beta}}d\xi, \quad t>0,$$
provided that the integral on the right-hand side of equality exists.

If $f\in L^p(0,T)$, $1\leq p\leq +\infty$, then (see Lemma 2.3 in \cite{Kilbas2006})
$$(I_{0t}^{\alpha}I_{0t}^{\beta} f ) (t)=(I_{0t}^{\alpha+\beta}f)(t).$$

For $0<\beta < 1$ an operator defined by  
\begin{equation}\label{Caputo}\partial_{0 t}^{\beta}g(t)=\frac{d}{dt}I_{0t}^{1-\beta}(g(t)-g(0))=\frac{1}{\Gamma(1-\beta)}\frac{d}{dt}\int\limits_{0}^{t}{\frac{g(\xi)-g(0)}{{{(t-\xi )}^{\beta}}}d\xi },\ \  t>0, \end{equation}		
is called the regularized Caputo fractional derivative, also known as the Dzherbashyan-Caputo fractional derivative \cite{{Dzherbashyan1968}}, provided that the right-hand side is well defined. Here $\Gamma(\beta)=\int\limits_{0}^{+\infty}x^{\beta-1}e^{-x}dx$ 
is  Euler's Gamma function. Here we suppose that the right-hand side of the equality \eqref{Caputo} exists.

If $g$ is absolutely continuous on $[0,T]$ one has  (\cite{Kilbas2006})
\begin{equation}\label{Caputo2}\partial_{0 t}^{\beta}g(t)=I_{0t}^{1-\beta}\frac{d}{dt}g(t)=\frac{1}{\Gamma(1-\beta)}\int\limits_{0}^{t}{\frac{g^\prime(\xi)}{{{(t-\xi )}^{\beta}}}d\xi },\ \  t>0. \nonumber \end{equation}	In addition, if $\partial_{0 t}^{\beta}g(t)\in C[0,T]$ then $g(t)-g(0)=I_{0 t}^{\beta}\partial_{0 t}^{\beta}g(t)$.

Next, we introduce the two-parameter function of the Mittag-Leffler type, which plays a very important role in the fractional calculus.
\subsection{The Mittag-Leffler function}

\begin{defn}(\cite{I.Podlubniy},\cite{Rudolf}) For $0<\alpha<1$, a two-parameter function of the Mittag-Leffler type is defined by the series expansion
\begin{equation}\label{ML-def}
{{E}_{\alpha ,\gamma  }(z)}=\sum\limits_{k=0}^{\infty }{\frac{{{z}^{k}}}{\Gamma (\alpha k+\gamma )}}, \quad		 z \in \mathbb{C}.
\end{equation}		
\end{defn}
If $\gamma= 1$, then the Mittag-Leffler function is called the one-parameter or classical Mittag-Leffler
function and it is denoted by $E_\alpha(z) = E_{\alpha,1}(z)$.

For $m=0,1,2,\ldots$, we denote by
$$
E_{\alpha,\gamma}^{(m)}(z)
:=\frac{d^m}{dz^m}E_{\alpha,\gamma}(z)
$$
the $m$th derivative of the Mittag--Leffler function with respect to its argument. By termwise differentiation of \eqref{ML-def}, one has (see \cite{Biolek2025})
\begin{equation}\label{ML-derivatives}
E_{\alpha,\gamma}^{(m)}(z)
=\sum_{j=0}^{\infty}
\frac{(j+m)!}{j!}
\frac{z^j}{\Gamma\bigl(\alpha(j+m)+\gamma\bigr)}.
\end{equation}
See also \cite{Rudolf,Gorenflo2002} for properties of the Mittag--Leffler function and its derivatives.

\begin{thm} 
  \label{ML-derivative-estimate}(see \cite[Theorem 1.6, p. 35]{I.Podlubniy} and \cite[Proposition 4, p. 34391; Proposition 6, p. 34393]{Biolek2025}).
Let $0<\alpha<2$, $\gamma\in\mathbb R$, and
$m=0,1,2,\ldots$. Let $\mu$ satisfy
$$
\frac{\pi\alpha}{2}<\mu<\min\{\pi,\pi\alpha\}.
$$
Then there exists a constant $C>0$ such that
$$
\left|E_{\alpha,\gamma}^{(m)}(z)\right|
\leq
\frac{C}{1+|z|^{m+1}},
\qquad
\mu\leq|\arg z|\leq\pi .
$$
\end{thm}

Now we present several  fundamental inequalities that are frequently employed in our analysis (\cite{Nirenberg}).

\begin{lem}\label{Gagliardo ineq}\cite{Nirenberg}
(Gagliardo-Nirenberg inequality).

Let $v\in L^q(\mathbb R^d)$ and $D^mv\in L^r(\mathbb R^d)$, where $1\leq q, r\leq \infty$. For $0\leq j<m$, the following inequality holds:
$$\|D^jv\|_{L^s}\leq C\|D^mv\|^\theta_{L^r}\|v\|^{1-\theta}_{L^q},$$
where $\frac{1}{s}=\frac{j}{d}+\theta(\frac{1}{r}-\frac{m}{d})+\frac{1-\theta}{q},$ $\frac{j}{m}\leq \theta.$

In the one dimensional case $d=1$, taking 
$j=0, m=1, q=r=2, s=2(p+1)$
we obtain $\theta=\frac{p}{2(p+1)}.$

Consequently,

$$\|v\|^{p+1}_{L^{2(p+1)}}\leq C\|v_x\|^{\frac{p}{2}}_{L^2}\|v\|^{\frac{p+2}{2}}_{L^2},$$
where $C>0$ is a constant independent of $v.$
\end{lem}
\begin{lem} \label{Lemm_yakupov}\cite{Yakupov1975}
If $u \in S(\mathbb R)$ and  $\int\limits_{-\infty}^{+\infty}\bigg(\frac{{\partial}^N u}{\partial x^N}\bigg)^2 dx\leq M<+\infty$, then the inequality is true.
\begin{equation}\label{eqlem7}
    \int\limits_{-\infty}^{+\infty}x^{2m}\bigg(\frac{\partial^k u}{\partial x^k}\bigg)^2 dx\leq c_1(k,m,M)\bigg(\int\limits_{-\infty}^{+\infty}x^{2m+2}u^2 dx\bigg)^{\frac{m}{m+1}}+c_2(k,m,M)\bigg(\int\limits_{-\infty}^{+\infty}x^{2m+2}u^2dx\bigg)^\frac{1}{2^k +1},\end{equation}
for $2mk\leq N$, where $c_j (k,m,M),  j=1,2$  are some positive numbers, while $m$, $k$ and $N$ are natural numbers.
\end{lem}

It is easy to see that from \eqref{eqlem7} it follows that
\begin{equation}\label{Asosiy lemma}   \int\limits_{-\infty}^{+\infty}x^{2m}\bigg(\frac{\partial^k u }{\partial x^k}\bigg)^2 dx\leq M_\varepsilon+\varepsilon \int\limits_{-\infty}^{+\infty}x^{2m+2}u^2 dx, \end{equation}
where $\varepsilon>0$ the positive constant $M_\varepsilon$ depends on $m,k,M$ and $N$.

Now we want to introduce semi-norms which are important in defining topological convergence in $S(\mathbb R)$.
According to  \cite{Yakupov1975} we introduce these semi-norms as follows
\begin{equation}\label{yarim norma}
|||u|||_{m,s}^2=\int\limits_{-\infty}^{+\infty}\bigg|\frac{\partial^m}{\partial x^m}u(x)\bigg|^2 dx+\int\limits_{-\infty}^{+\infty}\left(1+x^2\right)^s| u(x)|^2dx,
\end{equation}
where $m$ and $s$ are nonnegative integers.

Now we want to justify the above semi-norms. It is clear that if $v(x)\in S(\mathbb R)$ then the semi-norm $|||u|||_{m,s}$ is bounded for any $s,m\in \mathbb N\cup\{0\}.$ 

On the other hand, if $|||u|||_{m,s}$ is bounded for any $s,m\in \mathbb N\cup\{0\}$, then using the Cauchy–Bunyakovsky–Schwarz inequality and inequality \eqref{Asosiy lemma}, we have
$$ \sup_{x\in \mathbb R} \left|x^s\frac{d^m v(x)}{d x^m}\right|^2= 
2\sup_{x\in \mathbb R} \int\limits_{-\infty}^x\left( s\xi^{2s-1}\left(\frac{d^m v(\xi)}{d \xi^m}\right)^2
+\xi^{2s}\frac{d^{m} v(\xi)}{d \xi^{m}}\frac{d^{m+1} v(\xi)}{d \xi^{m+1}}\right)d\xi$$

$$\leq 2s\int\limits_{-\infty}^{+\infty} (1+x^2)^{s}\left(\frac{d^m v(x)}{d x^m}\right)^2dx$$

\begin{equation} \label{SN-SN}
    + 2\left(\int\limits_{-\infty}^{+\infty} (1+x^2)^{s}\left(\frac{d^m v(x)}{d x^m}\right)^2dx\int\limits_{-\infty}^{+\infty} (1+x^2)^{s}\left(\frac{d^{m+1} v(x)}{d x^{m+1}}\right)^2dx\right)^{1/2}\leq const.
\end{equation}

So, we can conclude that the condition \eqref{S(R)} in the definition of the Schwartz space $S(\mathbb R)$ is equivalent to the boundedness of the semi-norms $|||u|||_{m,s}$ for all $m, s\in \mathbb  N\cup{0}.$ 

We introduce the space of continuous functions $\mathbf u(t)=u(\cdot,t)$ as $SC(T)=C([0,T], S(\mathbb R))$. Following to \cite{Yakupov1975} we introduce the topological space $SC_\beta(T)=\{ u\in SC(T): \partial_{0t}^\beta u\in SC(T)\}$. Convergence in this space is defined by countably many semi-norms $\max_{0\leq t\leq T}|||u (\cdot,t)|||_{m,s}$, $\max_{0\leq t\leq T}|||(\partial_{0t}^\beta u) (\cdot,t)|||_{m,s},$ where $m, s =0,1,2,...\ .$
The following lemma is a direct consequence of Theorem 2.4 in \cite{Zacher2012}.
\begin{lem} \label{LemmAli}
\cite{Zacher2012}. Let $H$ be a real Hilbert space, $0<\beta<1$ and $T>0$. Suppose that 
$$v\in L^2([0,T];H), v-v_0=I_{0t}^\beta w \mbox{  for some  } w\in L^2([0,T];H) \mbox{ and } v_0\in H. $$
Then 
$$ (v(t),\partial_{0t}^\beta v(t))_H\geq\frac{1}{2}\partial_{0t}^\beta \|v(t)\|_H^2$$
for almost every $t\in [0,T].$
\end{lem}

The next statement follows from the comparison principle of Li and Liu \cite[Theorem 4.10]{LiLiu2018}, together with the explicit solution formula for the corresponding linear fractional equation \cite[Proposition 4.3]{LiLiu2018}.
\begin{lem}[Fractional Gronwall inequality]\label{GRONWALL}\cite{LiLiu2018}
Let $0<\beta<1$, $c\geq0$, and let $y\in C[0,T]$ and $q\in C[0,T]$ be nonnegative and $\partial_{0t}^\beta y\in C[0,T]$. Suppose that
\begin{equation}\label{weakGronwallAss}
\partial_{0t}^{\beta}y(t)\leq cy(t)+q(t)
\end{equation}
holds on $[0,T]$. Then
\begin{equation}\label{weakGronwallConclusion}
y(t)\leq y(0)E_\beta(ct^\beta)
+\int_0^t(t-s)^{\beta-1}E_{\beta,\beta}\!\left(c(t-s)^\beta\right)q(s)\,ds.
\end{equation}
Consequently, for $0\leq t\leq T$,
\begin{equation}\label{weakGronwallSimple}
y(t)\leq y(0)E_\beta(ct^\beta)
+\Gamma(\beta)E_{\beta,\beta}(ct^\beta)I_{0t}^{\beta}q(t).
\end{equation}
\end{lem}

\begin{lem}\label{yangi lemma1}[\cite{Ashurov};cf.\ \cite{YeGaoDing2007, KubicaYamamoto2018}]
Let $g_i(t),i=0,1,2,...,n$, $t\geq 0$ be a sequence of continuous  functions satisfying the inequality
$$g_i(t)\leq a+bI_{0t}^\beta g_{i-1}({t}), \quad i=1,2,...,n, $$
where $a\geq 0,b\geq 0, 0<\beta\leq 1,$ and $I_{0t}^\beta$ denotes the Riemann-Liouville fractional integral of order $\beta$. Then for all $n\in \mathbb{N}$, the following estimate holds
\begin{equation}\label{yangi lemma}
g_n (t)\leq a\sum_{i=0}^{n-1} \frac{b^i t^{i\beta}}{\Gamma(i\beta +1)}+b^n I_{0t}^{n\beta} g_0 (t).
\end{equation}
\end{lem} 
To investigate the solvability of the inverse problem for the nonlinear equation, we first consider the corresponding inverse problem for the linear equation. The solvability and stability results for this linear problem will be used in the construction of the iterative sequence for the nonlinear problem.

\section{INVERSE PROBLEM FOR THE LINEAR EQUATION}
We consider the following  inverse problem of finding the pair of functions $(u,f)$: 
\begin{align}    
Lu(x,t):=\partial _{0t}^\beta u(x,t)+a^2u_{xxxx}(x,t)+bu_{xxx}(x,t)+cu_{xx}(x,t)+du_x(x,t)+qu(x,t)=\nonumber \\
\label{u1,linear} =f(t)g(x,t)+h(x,t) \end{align}
with initial conditions
\begin{equation}
u(x,0)=\varphi(x), \quad x\in {\mathbb R}, \end{equation}
and an additional condition given by
\begin{equation}\label{u1,over.cond} 
\int\limits_{-     \infty}^{+\infty} \eta(x)u(x,t)dx=E(t),\quad 0\leq t\leq T,\end{equation}

$$\int\limits_{-\infty}^{+\infty}\eta(x)\varphi(x)dx=E(0),
$$

where $g,h\in SC(T)$, $E(t)\in C[0,T]$, $\varphi (x)\in S(\mathbb R)$.

Let $w$ is a solution to Cauchy problem $Lw=h, \ w(x,0)=\varphi(x)$.
Then can reduce our problem to the following inverse problem for $(v,f)$:
\begin{equation}\label{i,linear}
\partial _{0t}^\beta v(x,t)+a^2v_{xxxx}(x,t)+bv_{xxx}(x,t)+cv_{xx}(x,t)+dv_x(x,t)+qv(x,t)=f(t)g(x,t) \end{equation}
with initial conditions
\begin{equation}
v(x,0)=0, \end{equation}
and an additional condition given by
\begin{equation}\label{i,over.cond} 
\int\limits_{-\infty}^{+\infty} \eta(x)v(x,t)dx=\psi(t):=E(t)-\int\limits_{-\infty}^{+\infty}\eta(x)w(x,t)dx,\quad 0\leq t\leq T.\end{equation}
We multiply both sides of equation \eqref{i,linear} by the function $\eta(x)$ and integrate over $\mathbb R$, we obtain  
$$\partial _{0t}^\beta \int\limits_{-\infty}^{+\infty} \eta(x)v(x,t)dx+a^2\int\limits_{-\infty}^{+\infty} \eta(x)v_{xxxx}(x,t)dx+b\int\limits_{-\infty}^{+\infty}\eta(x)v_{xxx}(x,t)dx$$
\begin{equation}+c\int\limits_{-\infty}^{+\infty}\eta(x)v_{xx}(x,t)dx+d\int\limits_{-\infty}^{+\infty}\eta(x)v_x(x,t)dx+q\int\limits_{-\infty}^{+\infty}\eta(x)v(x,t)dx=f(t)\int\limits_{-\infty}^{+\infty}\eta(x)g(x,t)dx.  \end{equation}
We denote
$$g^* (t)=\int\limits_{-\infty}^{+\infty}\eta(x)g(x,t)dx,\quad \min_{t} |g^*(t)|=m>0.$$
After introducing the notation, we have the following second kind $f(t)$ operator equation 
\begin{equation}\label{i,f(t)}
f(t)=(Bf)(t)+\frac{\partial_{0t}^\beta \psi(t)}{g^*(t)},\end{equation}
where $\psi(t)$ is defined in \eqref{i,over.cond} and

$$(Bf)(t)=\frac{1}{g^*(t)}\bigg [a^2\int\limits_{-\infty}^{+\infty} \eta(x)v_{xxxx}dx+b\int\limits_{-\infty}^{+\infty}\eta(x)v_{xxx}dx$$
\begin{equation} \label{i,Bf}   
+c\int\limits_{-\infty}^{+\infty}\eta(x)v_{xx}dx+d\int\limits_{-\infty}^{+\infty}\eta(x)v_xdx+q\int\limits_{-\infty}^{+\infty}\eta(x)vdx\bigg].\end{equation}

\begin{thm}\label{i,thm}
    Let  $ \eta\in L^2(\mathbb R),$ and bounded, $\psi$, $\partial^\beta_{0t}\psi\in C[0,T]$, $h, g\in SC(T)$, $g^*(t)>0$,    then the inverse problem  \eqref{u1,linear}-\eqref{u1,over.cond} has a unique solution $(u,f)\in SC_\beta(T)\times C[0,T]$ and the following stability estimate holds 
     $$\|f\|_{C{[0,T]}}\leq M \|\partial ^\beta_{0t}\psi\|_{C{[0,T]}},$$
     where $M$ is a positive constant.
\end{thm}
\begin{proof}
We prove that equation \eqref{i,f(t)} admits a solution. We need to study the operator $B.$ The operator $B: C[0,T] \rightarrow C[0,T].$ 

 Now we need to obtain some a-priori estimates. For this, first we prove      the following lemma.

\begin{lem}
    Let $T>0$  and suppose that the conditions of Theorem \ref{i,thm} hold. Then 
   
    $$\left\| \frac{\partial^n}{\partial x^n}v(\cdot,t)\right\|^2\leq C \Gamma(\beta) E_{\beta,\beta}(C_0t^\beta)I^{\beta}_{0t} |f(t)|^2,$$
    
    for all $0\leq t\leq T$, where $n$ is a non negative integer, $C$ is a positive constant, and 
    $C_0=C_0(a,c,q)$. 
  \end{lem}
Let us now carry out an induction on the order of the derivative. We put 
$$|||v|||_n^2=\int\limits_{-\infty}^{+\infty} \left( \frac{d^n v(x)}{d x^n}\right)^2dx,$$
and $|||v|||_0=\|v\|_{L^2(\mathbb R)}$.

Now we can take the $n$th derivative of each term of the equation with respect to $x$, then multiply by $2\frac{\partial^nv}{\partial x^n}$, and integrate with respect to $x$ in $\mathbb{R}$.
\begin{equation}\label{i,equation}
2\int\limits_{-\infty}^{+\infty}\partial_{0t}^\beta \frac{\partial^n v}{\partial x^n}\frac{\partial^n v}{\partial x^n}dx+2\int\limits_{-\infty}^{+\infty}\frac{\partial^{n} }{\partial x^{n}}Lv{\frac{\partial^n v}{\partial x^n }}dx
=2f(t)\int\limits_{-\infty}^{+\infty}\frac{\partial^n g}{\partial x^n}\frac{\partial^n v}{\partial x^n}dx,\end{equation}
where $Lv=a^2v_{xxxx}+bv_{xxx}+cv_{xx}+dv_x+qv.$

Now we estimate the terms on the left-hand side of this equality.
\begin{equation}\label{i,partial}
2\int\limits_{-\infty}^{+\infty}\frac{\partial^n v}{\partial x^n}\partial_{0t}^\beta \frac{\partial^n v}{\partial x^n}dx\geq\partial_{0t}^\beta
\int\limits_{-\infty}^{+\infty}\left(\frac{\partial^n v}{\partial x^n}\right)^2 dx=\partial_{0t}^\beta|||v|||_n^2, \end{equation}
\begin{equation}\label{i,xxxx}
2a^2\int\limits_{-\infty}^{+\infty}\frac{\partial^{n+4} v}{\partial x^{n+4}}\frac{\partial^n v}{\partial x^n}dx=2a^2\int\limits_{-\infty}^{+\infty}\left(\frac{\partial^{n+2} v}{\partial x^{n+2}}\right)^2 dx=2a^2|||v|||^2_{n+2},\end{equation}
\begin{equation}\label{i,toq}
\int\limits_{-\infty}^{+\infty}\frac{\partial^{n+3}v}{\partial x^{n+3}}\frac{\partial^n v}{\partial x^n}dx=0, \ \ \int\limits_{-\infty}^{+\infty}\frac{\partial^{n+1}v}{\partial x^{n+1}}\frac{\partial^n v}{\partial x^n}dx=0, \end{equation}

\begin{equation}\label{i,xx}
\left|2\int\limits_{-\infty}^{+\infty}\frac{\partial^{n+2}v}{\partial x^{n+2}}\frac{\partial^n v}{\partial x^n}dx\right|=2|||v|||_{n+1}^2\leq \sigma |||v|||_{n+2}^2+ \frac{1}{\sigma}|||v|||_n^2. \end{equation}
The last integral on the right-hand side of \eqref{i,equation} is estimated using the Cauchy inequality:
\begin{equation}\label{i,f}
\left|2\int\limits_{-\infty}^{+\infty}f(t)\frac{\partial^n g}{\partial x^n}\frac{\partial^n v}{\partial x^n}dx\right|
\leq (f(t))^2\int\limits_{-\infty}^{+\infty}\left(\frac{\partial^n g}{\partial x^n}\right)^2 dx+\int\limits_{-\infty}^{+\infty} \left(\frac{\partial^n v}{\partial x^n}\right)^2 dx.   \end{equation}
Summarizing the above estimate, using the appropriate choice of $\sigma=\frac{a^2}{2}$, we have
\begin{equation}\label{i,gronwall}
\partial_{0t}^\beta |||v|||_n^2 \leq C_0|||v|||_n^2+C|f(t)|^2.\end{equation}
Using  Gronwall's inequality (see Lemma\ref{GRONWALL})  we obtain the desired estimate
\begin{equation}\label{i,gronwall yechimi}
|||v|||_n^2\leq C\Gamma(\beta) E_{\beta,\beta}(C_0t^\beta)I^{\beta}_{0t} |f(t)|^2. \end{equation}
Now we study the properties of the operator $B$.
For that, we need the following estimates.
$$I_k=\int\limits_{-\infty}^{+\infty}\eta(x)\frac{\partial^kv{(x,t)}}{\partial x}dx, \quad k=0,1,2,3,4. $$
Using Cauchy's inequality, we obtain the following estimate for each $I_k$
\begin{equation}\label{i,I_k}
|I_k|\leq \|\eta\|_{L^2}\cdot \bigg \|\frac{\partial^k v}{\partial x}\bigg \|_{L^2}\leq C(T)\sqrt{I^\beta_{0t}|f(t)|^2}. \end{equation} 
Using the definition of $Bf$ and expression \eqref{i,I_k}, we obtain the following
$$|Bf|^2\leq \frac{5}{m} \bigg(\sum^4_{k=0}I^2_k\bigg)\leq \tilde C 
I^\beta_{0t}|f(t)|^2. $$
Second iteration 
$$|B^2f|^2=|B(B(f))|^2\leq \tilde CI^\beta_{0t}|Bf|^2\leq \tilde  C I^\beta_{0t}\circ (\tilde C I^\beta _{0t}|f(t)|^2)\leq \tilde C^2I^{2\beta}_{0t}|f(t)|^2.$$
Continuing in the same manner, we repeat this process $n$ times. From this we find the norm of the operator $B^n.$
\begin{equation}\label{i,B^n}
|B^n f|^2=|B(B^{n-1}f)|^2\leq \tilde C^n I^{n\beta}_{0t}|f(t)|^2.\end{equation}
Next, we estimate the norm of $Bf$.
\begin{equation}
\|Bf\|_{C[0,T]}^2\leq \tilde C I^\beta_{0t}|f|^2\leq \tilde C \max _{0\leq t\leq T}|f|^2\cdot I^\beta_{0t}1\leq\frac{\tilde CT^{\beta}}{\Gamma(\beta+1)} \|f\|^2_{C[0,T]}. \end{equation}
Now, iterating the inequality \ref{i,B^n}
$$\|B^nf\|^2_{C[0,T]}\leq\frac{\tilde C^nT^{n\beta}}{\Gamma(n\beta+1)}\cdot \|f\|^2_{C[0,T]}.$$
By definition
$$\|B^n\|=\sup_{\|f\|=1}\|B^n f\|_{C[0,T]}\leq \sqrt{\frac{\tilde C^n T^{n\beta}}{\Gamma(n\beta+1)}}  .$$
According to Stirling's formula,
$$\sqrt{\frac{\tilde C^n T^{n\beta}}{\Gamma(n\beta+1)}}\sim \sqrt{\frac{(\tilde CT^\beta)^n}{\sqrt{2\pi n\beta}\cdot (\frac{n\beta}{e})^{n\beta}}}, \quad n\rightarrow{+\infty}.$$
 Now, let us calculate the spectral radius
$$r(B)=\lim_{n \rightarrow {+\infty}}\sqrt[n]{\|B^n\|}\leq \lim_{n \rightarrow {+\infty}}\frac{\tilde CT^{\frac{\beta}{2}}}{(2\pi \beta)^{\frac{1}{4n}}\sqrt[4]{\sqrt[n]{n}}(\frac{n\beta}{e})^{{\frac{\beta}{2}}}}=0.$$
Therefore, $B$ is a quasinilpotent operator \cite{Kreyszig}. Consequently, the spectrum of the operator $B$ consists only of the point $\{\lambda=0\}$, so it follows that the resolvent operator $(I-B)^{-1}$ exists. 
Hence, 
$$\|(I-B)^{-1}f\|_{C[0,T]}\leq \sum _{k=0}^{+\infty}\sqrt{\frac{\tilde C^k T^{k\beta}}{\Gamma(k\beta+1)}}  \|f\|_{C[0,T]}.$$
So, summarizing the above, the resolvent operator 
$$(I-B)^{-1}:C[0,T]\rightarrow C[0,T]$$
is bounded and continuous. From this, according to \eqref{i,f(t)},
 $f$ can be determined
$$f=(I-B)^{-1}\frac{\partial_{0t}^\beta \psi}{g^*}.$$
Then, according to the estimate for the resolvent operator, we have 
 $$\|f\|_{C[0,T]}\leq M \|\partial ^\beta_{0t}\psi\|_{C{[0,T]}},$$
where $M=\sum _{k=0}^{+\infty}\sqrt{\frac{\tilde C^k T^{k\beta}}{\Gamma(k\beta+1)}}<{+\infty},$\quad since $r(B)=0.$

We proved the existence of unique $f\in C[0,T].$ The existence of unique $v\in SC_\beta(T)$ follows from  (see  Theorem 3.1. in \cite{Ashurov}). Hence, the theorem is proved.

It has been established that $w\in SC_\beta(T)$ is the solution of the linear direct problem and $v\in SC_\beta (T)$ is the solution of the linear inverse problem. Since $u=w+v$, it follows that $u\in SC_\beta (T)$. Therefore, $(u,f)\in SC_\beta(T)\times C[0,T].$ 
\end{proof}

\section{SOLVABILITY OF THE INVERSE PROBLEM FOR THE NONLINEAR EQUATION} 
Now we consider the inverse problem to find $(u,f)$ for the  Kuramoto-Sivashinsky equation 
\begin{equation}\label{41}
\partial _{0t}^{\beta }u+a^2u_{xxxx}+bu_{xxx}+cu_{xx}+du_x+qu=f(t)g(x,t)+\gamma u^pu_x, \quad x\in {\mathbb R},\quad 0<   t\le T,\end{equation}		
with initial condition	
\begin{equation}\label{41-initial}u(x,0)=\varphi(x), \quad x\in {\mathbb R}, \end{equation} 
and the additional integral condition
\begin{equation}\label{41-addC}
\int\limits_{-\infty}^{+\infty}\eta(x)u(x,t)dx=E(t),\quad 0\leq   t\leq T,
\end{equation}

$$\int\limits_{-\infty}^{+\infty}\eta(x)\varphi(x)dx=E(0),
$$
where  $0<\beta <1$, $a>0$, $b, c, d, q$ and $\gamma$ are real numbers, \quad $p\in\mathbb N$.

In this paragraph, we prove the following local solvability theorem.

\begin{thm}\label{th4.1}  Let $\eta\in W_2^4(\mathbb{R})$,
$E,\partial_{0t}^{\beta}E\in C[0,T]$, $g\in SC(T)$, and
\[
g^*(t)\geq m>0,\qquad 0\leq t\leq T.
\]
Then there exists a number $0<t_1\leq T$, dependent on the
coefficients of the equation, $\|\eta\|_{W_2^4(\mathbb{R})}$, $\|\varphi\|_1$,
$\|\partial_{0t}^{\beta}E\|_{C[0,T]}$, $m^{-1}$, and $\|g\|_{C([0,T];L^2(\mathbb{R}))},$
such that the inverse problem \eqref{41}--\eqref{41-addC} has a unique solution
$(u,f)\in SC_{\beta}(t_1)\times C[0,t_1].$
\end{thm} 
   
We solve this problem by successive approximation: $u_0(x,t)=\varphi(x),$ 
\begin{equation} \label{equation}
\partial _{0t}^\beta u_n+a^2u_{nxxxx}+bu_{nxxx}+cu_{nxx}+du_{nx}+qu_n=f_n(t)g(x,t)+\gamma u^p_{n-1}u_{n-1,x}, n\geq 1,
\end{equation}
with initial conditions
$$u_n(x,0)=\varphi(x),\quad x\in \mathbb R,$$
and an additional condition given by
\begin{equation} \label{test}
\int\limits_{-\infty}^{+\infty} \eta(x)u_n(x,t)dx=E(t),\quad 0\leq t\leq T.
\end{equation}

According to Theorem \ref{i,thm} $u_n\in SC_\beta(T)$, $f_n\in C[0,T]$. To prove solvability of the problem \eqref{41}-\eqref{41-initial}, \eqref{41-addC} we need to prove the convergence of the sequence $\{(u_n, f_n)\}$. 

We multiply both sides of equation \eqref{equation} by the function $\eta(x)$ and integrate over $\mathbb R$, we obtain  
$$\partial _{0t}^\beta \int\limits_{-\infty}^{+\infty} \eta(x)u_ndx+a^2\int\limits_{-\infty}^{+\infty} \eta(x)u_{nxxxx}dx+b\int\limits_{-\infty}^{+\infty}\eta(x)u_{nxxx}dx$$
\begin{equation}+c\int\limits_{-\infty}^{+\infty}\eta(x)u_{nxx}dx+d\int\limits_{-\infty}^{+\infty}\eta(x)u_{nx}dx+q\int\limits_{-\infty}^{+\infty}\eta(x)u_ndx=f_n(t)\int\limits_{-\infty}^{+\infty}\eta(x)g(x,t)dx+\gamma \int\limits_{-\infty}^{+\infty}\eta(x)u^p_{n-1}u_{n-1,x}dx.
\end{equation}


We denote

\begin{equation}
g^* (t)=\int\limits_{-\infty}^{+\infty}\eta(x)g(x,t)dx,\end{equation}
after the notation, we obtain 
$$f_n(t)=(Bf_n)(t)+\frac{\partial_{0t}^\beta E(t)}{g^*(t)},$$
where $E(t)$ is defined in \eqref{test} and

$$(Bf_n)(t)=\frac{1}{g^*(t)}\bigg [a^2\int\limits_{-\infty}^{+\infty} \eta(x)u_{nxxxx}dx+b\int\limits_{-\infty}^{+\infty}\eta(x)u_{nxxx}dx$$
\begin{equation}   \label{B_n} 
+c\int\limits_{-\infty}^{+\infty}\eta(x)u_{nxx}dx+d\int\limits_{-\infty}^{+\infty}\eta(x)u_{nx}dx+q\int\limits_{-\infty}^{+\infty}\eta(x)u_ndx+\frac{\gamma}{p+1} \int\limits_{-\infty}^{+\infty}\eta_x(x)u^{p+1}_{n-1}dx
\bigg].
\end{equation}
The existence of a solution follows from Theorem \ref{i,thm} for each $n.$ Now we show that the sequence converges in the topology defined by the semi-norms \eqref{yarim norma}. So, we need to obtain countably many a priori estimates in the semi-norms \eqref{yarim norma}. Based on these estimates, at the end, we show that $\{u_n\}$ is a Cauchy sequence.

Now we establish some auxiliary lemmas.

We put 
$$\|v\|^2=\int\limits_{-\infty}^{+\infty}v^2(x)dx.$$

The following lemma directly follows from \eqref{B_n}.
\begin{lem}
    Let $g^*(t)>m>0$ and the conditions of Theorem \ref{th4.1} are valid. Then for $t_1\in (0,T)$ the following estimate holds 
    \begin{equation} \label{estimate f_n}
    |f_n(t)|\leq \frac{1}{m}\max_{0\leq t\leq t_1}|\partial^\beta_{0t}E(t)|+ \tilde  c_1\bigg(\|u_n\|^2+\|u_{nx}\|^2\bigg)^{1/2}+\tilde c_2\bigg(\|u_{n-1}\|^2+\|u_{n-1,x}\|^2\bigg)^{(p+1)/2},\end{equation}
where $0\leq t\leq t_1$, $ \tilde c_1,\tilde c_2>0$ are constants that depend only on the coefficients, $m$ and $\|\eta\|_{W^4_2(\mathbb R)}.$

\end{lem}
 






 \subsection{Estimate for $\|u_n\|^2+\|u_{nx}\|^2$.}
     
\begin{lem} \label{norma2}
    There exists a number $0<t_1\leq T$ dependent on coefficients of the equation, $ \|\varphi\|_1$, $\|\partial^\beta_{0t}E\|_{C[0,T]}$ and $\min_t\|g\|_{L^2}$, such that
  \begin{equation}\label{U}  
  \max_{0\leq t 
\leq t_1}(\|u_n\|^2+\|u_{nx}\|^2)\leq const <+\infty,\end{equation}
\begin{equation} \label{max f_n}
\max_{0\leq t 
\leq t_1}|f_n(t)|\leq const<+\infty.\end{equation}
           where the constant does not depend on $n$.
\end{lem}
\begin{proof}    
Further we use $C$ for different constants, dependent on $m$, $\|\eta\|_{W^4_2(\mathbb R)}$ and on the coefficients of the equation.

Multiply equation \eqref{equation} by $2u_n$ and integrate with respect to $x$ over ${\mathbb R}$.

$$2 \int\limits_{-\infty}^{+\infty}u_n(x,t)\partial _{0t}^\beta u_n (x,t)dx+2a^2\int\limits_{-\infty}^{+\infty} u_n (x,t)u_{nxxxx}(x,t)dx+2b\int\limits_{-\infty}^{+\infty} u_n(x,t)u_{nxxx}(x,t)dx$$
$$+2c\int\limits_{-\infty}^{+\infty}u_n(x,t)u_{nxx}(x,t)dx+2d\int\limits_{-\infty}^{+\infty}u_n(x,t)u_{nx}(x,t)dx+2q\int\limits_{-\infty}^{+\infty}u_n(x,t)u_n(x,t)dx$$
\begin{equation} \label{u_{n}}
=2f_n(t)\int\limits_{-\infty}^{+\infty}u_n(x,t)g(x,t)dx-\frac{2\gamma}{p+1} \int\limits_{-\infty}^{+\infty}u_{nx}(x,t)u^{p+1}_{n-1}(x,t)dx.
\end{equation}
We now estimate the terms on the left side of \eqref{u_{n}}.

$$2\int\limits_{-\infty}^{+\infty}u_n \partial^\beta _{0t}u_ndx\geq  \partial^\beta_{0t}\|u_n\|^2,\quad 2a^2\int\limits_{-\infty}^{+\infty}u_nu_{nxxxx}dx=2a^2\|u_{nxx}\|^2,$$
$$\int\limits_{-\infty}^{+\infty}u_nu_{nxxx}dx=0,\quad |\int\limits_{-\infty}^{+\infty}u_nu_{nxx}dx|=\|u_{nx}\|^2,$$

$$\int\limits_{-\infty}^{+\infty}u_nu_{nx}dx=0,\quad\int\limits_{-\infty}^{+\infty}u^2_ndx=\|u_n\|^2.$$
We estimate the first term on the right-hand side of the equation. 

$$2\bigg|\int\limits_{-\infty}^{+\infty}u_n(x,t)f_n(t)g(x,t)dx\bigg|\leq 2\|u_n\||f_n|\|g\| \leq \|u_n\|^2+|f_n(t)|^2\|g\|^2.$$

Now we examine the second term on the right-hand side of the equation  

$$\bigg|2\gamma\int\limits_{-\infty}^{+\infty} u_nu^p_{n-1}u_{n-1,x}dx\bigg|=\bigg|-\frac{2\gamma}{p+1}\int\limits_{-\infty}^{+\infty}u_{nx}u^{p+1}_{n-1}dx\bigg|\leq \frac{2|\gamma|}{p+1}\|u_{nx}\|\|u_{n-1}\|_{L^{2(p+1)}}^{p+1}$$
$$\leq \frac{C}{p+1}\|u_{nx}\|\|u_{n-1}\|^{\frac{p+2}{2}}\|u_{n-1,x}\|^{\frac{p}{2}}\leq \|u_{nx}\|^2+C\|u_{n-1}\|^{p+2}\|u_{n-1,x}\|^p$$
$$\leq \|u_{nx}\|^2+C\bigg(\|u_{n-1}\|^2+\|u_{n-1,x}\|^2\bigg)^{p+1}.$$

Summarizing the above estimates,

$$\partial^\beta_{0t}\|u_n\|^2+2a^2\|u_{nxx}\|^2\leq (2|c|+1)\|u_{nx}\|^2+(2|q|+1)\|u_n\|^2$$
\begin{equation} \label{partial u_n}
+|f_n|^2\|g\|^2 +C\bigg(\|u_{n-1}\|^2+\|u_{n-1,x}\|^2\bigg)^{p+1}.\end{equation}

Now, we take the first order $x$-derivative and multiply both sides of equation \eqref{equation} by $2u_{nx}$: 
$$2\partial _{0t}^\beta \int\limits_{-\infty}^{+\infty} u_{nx} u_{nx}dx+2a^2\int\limits_{-\infty}^{+\infty} u_{nx}u_{nxxxxx}dx+2b\int\limits_{-\infty}^{+\infty} u_{nx}u_{nxxxx}dx$$
$$+2c\int\limits_{-\infty}^{+\infty}u_{nx}u_{nxxx}dx+2d\int\limits_{-\infty}^{+\infty}u_{nx}u_{nxx}dx+2q\int\limits_{-\infty}^{+\infty}u_{nx}u_{nx}dx$$
\begin{equation} \label{equation u_{nx}}
=2f_n(t)\int\limits_{-\infty}^{+\infty}u_{nx}g_xdx+2\gamma \int\limits_{-\infty}^{+\infty}(u^p_{n-1}u_{n-1,x})_xu_{nx}dx.
\end{equation}

Here we use similar inequalities as in the previous case.
       
$$2\int\limits_{-\infty}^{+\infty}u_{nx} \partial^\beta _{0t}u_{nx}dx\geq  \partial^\beta_{0t}\|u_{nx}\|^2,\quad 2a^2\int\limits_{-\infty}^{+\infty}u_{nx}u_{nxxxxx}dx=2a^2\|u_{nxxx}\|^2,$$
$$\int\limits_{-\infty}^{+\infty}u_{nx}u_{nxxxx}dx=0,\quad \int\limits_{-\infty}^{+\infty}u_{nx}u_{nxx}dx=0, $$

$$\bigg|\int\limits_{-\infty}^{+\infty}u_{nx}u_{nxxx}dx\bigg|=\|u_{nxx}\|^2\leq \|u_{nx}\|\|u_{nxxx}\|\leq \varepsilon \|u_{nxxx}\|^2+\frac{1}{4\varepsilon}\|u_{nx}\|^2.$$

We estimate the terms on the right-hand side of the equation

$$\bigg|\int\limits_{-\infty}^{+\infty}f_n(t)g_xu_{nx}dx\bigg|=\bigg|f_n(t)\int\limits_{-\infty}^{+\infty}gu_{nxx}dx\bigg|\leq \frac{1}{2}\|u_{nxx}\|^2+\frac{1}{2}|f_n|^2\|g\|^2$$
$$\leq \frac{1}{2}\|u_{nx}\|\|u_{nxxx}\|+\frac{1}{2}|f_n|^2\|g\|^2\leq \varepsilon \|u_{nxxx}\|^2+\frac{1}{8\varepsilon}\|u_{nx}\|^2+\frac{1}{2}|f_n|^2\|g\|^2.$$

$$\bigg|\frac{2\gamma}{p+1}\int\limits_{-\infty}^{+\infty}u_{nxxx}u^{p+1}_{n-1}dx\bigg|\leq \frac{2|\gamma|}{p+1}\|u_{nxxx}\|\|u_{n-1}\|^{p+1}_{L^{2(p+1)}}\leq{C}\|u_{nxxx}\|\|u_{n-1}\|^{\frac{p+2}{2}}\|u_{n-1,x}\|^{\frac{p}{2}}$$
$$\leq \varepsilon\|u_{nxxx}\|^2+C\|u_{n-1}\|^{p+2}\|u_{n-1,x}\|^p\leq \varepsilon\|u_{nxxx}\|^2+C\bigg(\|u_{n-1}\|^2+\|u_{n-1,x}\|^2\bigg)^{p+1}.$$

From the above estimates, it follows that
$$\partial^\beta_{0t}\|u_{nx}\|^2+2a^2\|u_{nxxx}\|^2\leq (2|c|+1)\|u_{nxx}\|^2+2|q|\|u_{nx}\|^2$$
\begin{equation} \label{partial u_{nx}}
 +|f_n|^2\|g\|^2 +\varepsilon\|u_{nxxx}\|^2+C\bigg(\|u_{n-1}\|^2+\|u_{n-1,x}\|^2\bigg)^{p+1}.\end{equation}

If we add expressions \eqref{partial u_n} and \eqref{partial u_{nx}}, the following
$$
\partial^\beta_{0t}\bigg\{\|u_n\|^2+\|u_{nx}\|^2\bigg\}+(2a^2-2(1+|c|)\varepsilon)\|u_{nxxx}\|^2+2a^2\|u_{nxx}\|^2$$
\begin{equation}\label{4.16}
\leq C\|u_{nx}\|^2+C\|u_n\|^2+C\bigg(\|u_{n-1}\|^2+\|u_{n-1,x}\|^2\bigg)^{p+1}+|f_n|^2\|g\|^2.\end{equation}
 Now from \eqref{estimate f_n}, \eqref{4.16} and choosing $\varepsilon=\frac{a^2}{2(1+|c|)}$ we obtain
\begin{equation}
\partial^\beta_{0t}\bigg\{\|u_n\|^2+\|u_{nx}\|^2\bigg\}\leq C_1\bigg\{\|u_n\|^2+\|u_{nx}\|^2\bigg\}+C_2\bigg\{\|u_{n-1}\|^2+\|u_{n-1,x}\|^2\bigg\}^{p+1}+\frac{1}{m^2}\max_{0\leq t\leq T}|\partial^\beta_{0t}E(t)|^2\|g\|^2.
\end{equation}

Now we denote
$Y_n=\|u_n\|^2+\|u_{nx}\|^2,\quad Y_0=\|\varphi\|^2+\|\varphi_{x}\|^2.$ Then

$$\partial^\beta_{0t}Y_n\leq C_1Y_n+C_2Y^{p+1}_{n-1}+\frac{1}{m^2}\max_{0\leq t\leq T}|\partial^\beta_{0t}E(t)|^2\|g\|^2,\quad n\geq 1.$$
Using Gronwall's inequality from Lemma \ref{GRONWALL}, we have
$$Y_n(t)\leq Y_0 E_\beta(C_1t^\beta)+\Gamma(\beta)E_{\beta,\beta}(C_1t^\beta)I^\beta_{0t}\bigg(C_2Y^{p+1}_{n-1}+\frac{1}{m^2}\max_{0\leq t\leq T}|\partial^\beta_{0t}E(t)|^2\|g\|^2\bigg). $$

We put
$$A_0=Y_0E_\beta(C_1T^\beta)+\frac{T^\beta}{\beta m^2}E_{\beta,\beta}(C_1T^\beta)\max_{0\leq t\leq T}|\partial^\beta_{0t}E(t)|^2\|g\|^2.$$
It is clear $Y_0\leq 2A_0.$
Let $$Y_{n-1}(t)\leq 2A_0, \quad n\geq1.$$
Then
\begin{equation}
Y_n\leq \frac{2^{p+1}C_2t^\beta}{\beta}E_{\beta,\beta}(C_1T^\beta)A^{p+1}_0+A_0.
\end{equation}
Choosing 
$$t_1\leq \bigg(\frac{\beta}{2^{p+1}A^p_0C_2E_{\beta,\beta}(C_1T^\beta)}\bigg)^{\frac{1}{\beta}},$$
we get
$\max_{0\leq t\leq t_1}Y_n\leq 2A_0.$
Therefore, inequality $\|u_n\|^2+\|u_{nx}\|^2\leq 2A_0 $  holds for all $n$ and $0<t<t_1\leq T.$

Estimate \eqref{max f_n} follows from \eqref{estimate f_n} and \eqref{U}.

\end{proof}


\subsection{Estimates for higher derivatives.}
\begin{lem} \label{high derivatives}
    Let $t_1$ be the same as in Lemma \ref{norma2}. Then 
    $$\left\| \frac{\partial^k}{\partial x^k}u_n(\cdot,t)\right\|^2\leq M(k),$$
    for all $0\leq t\leq t_1$, where $k\geq 2$ and $M(k)$ is a positive constant.
\end{lem}
\begin{proof}
    
Let us now carry out an induction on the order of the derivative. For $k\geq 1$ we put 
$$|||v|||_k^2=\int\limits_{-\infty}^{+\infty} \left( \frac{d^k v(x)}{d x^k}\right)^2dx,$$
and $|||v|||_0=\|v\|_{L^2}$.

We suppose $|||u_n|||_j\leq M(j)<+\infty$ for all $0\leq j<k$, $n\in \mathbb N$, where $k\geq 2$.  Now we can take the $k$th derivative of each term of the equation with respect to $x$, then multiply by $2\frac{\partial^k u_n}{\partial x^k}$, and integrate with respect to $x$ in $\mathbb{R}$.

$$2\int\limits_{-\infty}^{+\infty}\partial_{0t}^\beta \frac{\partial^k u_n}{\partial x^k}\frac{\partial^k u_n}{\partial x^k}dx+2\int\limits_{-\infty}^{+\infty}\frac{\partial^{k} }{\partial x^{k}}Lu_n {\frac{\partial^k u_n}{\partial x^k }}dx$$

\begin{equation}\label{LEYBNITS}
=2\gamma{\int\limits_{-\infty}^{+\infty} {\frac{\partial^k}{\partial x^k}}
{\left(u^p_{n-1}\frac{\partial {u_{n-1}}}{\partial x}\right)}\frac{\partial^k u_n}{\partial x^k}dx 
+2f_n(t)\int\limits_{-\infty}^{+\infty}\frac{\partial^k g}{\partial x^k}\frac{\partial^k u_n}{\partial x^k}dx}.\end{equation}






We estimate only the integral arising from the nonlinear term. Other terms in \eqref{LEYBNITS} can be estimated as in the proof of Lemma \ref{norma2}.
\begin{equation}\label{4.24}
\int\limits_{-\infty}^{+\infty}\frac{\partial^k}{\partial x^k} \left(u^p_{n-1} \frac{\partial u_{n-1}}{\partial x}\right)\frac{\partial^k u_n}{\partial x^k}dx
=\frac{1}{p+1}\int\limits_{-\infty}^{+\infty}\frac{\partial^k u_n}{\partial x^k}\frac{\partial^{k+1} }{\partial x^{k+1}}(u^{p+1}_{n-1})dx
=\frac{1}{p+1}\int\limits_{-\infty}^{+\infty}\frac{\partial^{k+2}u_n}{\partial x^{k+2}} \frac{\partial^{k-1}} {\partial x ^{k-1}}{(u^{p+1}_{n-1})}dx.\end{equation}

If we apply Leibniz's formula, we get the following

\begin{equation}
\frac{\partial^{k-1}}{\partial x^{k-1}}(u^{p+1}_{n-1})=\sum _{\alpha_1+\alpha_2+...+\alpha_{p+1}=k-1}\frac{(k-1)!}{\alpha_1!...\alpha_{p+1}!}\prod^{p+1}_{i=1}\frac{\partial^{\alpha_i}u_{n-1}}{\partial x^{\alpha_i}}.\end{equation}
Let us define

$$I=\int\limits_{-\infty}^{+\infty}\frac{\partial^{k+2}}{\partial x^{k+2}}u_n\frac{\partial^{k-1}}{\partial x^{k-1}}(u^{p+1}_{n-1})dx$$

\begin{equation}
=\sum _{\alpha_1+\alpha_2+...+\alpha_{p+1}=k-1}\frac{(k-1)!}{\alpha_1!...\alpha_{p+1}!}\int\limits_{-\infty}^{+\infty}\frac{\partial^{k+2}u_n}{\partial x^{k+2}}\prod^{p+1}_{i=1}\frac{\partial^{\alpha_i}u_{n-1}}{\partial x^{\alpha_i}}dx.\end{equation}

It is possible the following two cases.

\textbf{First case:}  $\alpha_i\leq k-2,$ \quad for all $i.$ Then
\begin{equation} 
\sup_{x\in \mathbb R}\left|\frac{\partial^{\alpha_j} u_{n-1}}{\partial x^{\alpha_j}}\right|\leq\sqrt{2}\left(\int\limits_{-\infty}^{+\infty}\left(\frac{\partial^{\alpha_j} u_{n-1}}{\partial x^{\alpha_j}}\right)^2dx\int\limits_{-\infty}^{+\infty}\left(\frac{\partial^{\alpha_j+1}u_{n-1}}{\partial x^{\alpha_j+1}}\right)^2 dx\right)^{\frac{1}{4}}\leq C.\end{equation}

By induction hypotesis $|||u_{n-1}|||_{\alpha_i}\leq M(\alpha_i)$. So, 
$$|I|\leq C |||u_n|||_{k+2}|||u_{n-1}|||_{\alpha_i}\leq \varepsilon|||u_n|||_{k+2}^2 +C(\varepsilon).$$

\textbf{Second case:} $\alpha_j=k-1$ for some $1\leq j\leq p+1$ and 
$\alpha_i=0$, \quad $i\neq j.$

In this case, using Cauchy-Schwarz, Gagliardo-Nirenberg \ref{Gagliardo ineq}, and H\'older inequalities, we have
$$|I|=\bigg|\int\limits_{-\infty}^{+\infty}\frac{\partial ^{k+2}u_n}{\partial x^{k+2}}u^p_{n-1}\frac{\partial^{k-1}}{\partial x^{k-1}}u_{n-1}dx\bigg|\leq \sup_{x\in {\mathbb R}}\bigg|\frac{\partial^{k-1} u_{n-1}}{\partial x^{k-1}}\bigg|\bigg\|\frac{\partial ^{k+2}u_n}{\partial x^{k+2}}\bigg\|_{L^2}\|u_{n-1}\|^p_{L^{2p}}$$
\begin{equation}\label{h derivative}
\leq C|||u_n|||_{k+2}(1+|||u_{n-1}|||_k)\leq \varepsilon |||u_n|||_{k+2}^2+C(\varepsilon)(1+|||u_{n-1}|||_k^2),\end{equation}
where $\varepsilon>0$ and its appropriate value will be chosen later.



The last integral on the right-hand side of \eqref{LEYBNITS} is estimated using the Cauchy inequality:

\begin{equation}
\left|2f_n(t)\int\limits_{-\infty}^{+\infty} \frac{\partial^k g}{\partial x^k}\frac{\partial^k u_n}{\partial x^k}dx\right|
\leq M||| g|||_k^2 +M ||| u_n|||_k^2 .\end{equation}            

Summarizing the above estimation, and choosing of $\varepsilon$ small then $a^2$, we have
\begin{equation}
\partial_{0t}^\beta |||u_n|||_k^2 \leq C( |||u_{n}|||_k^2+|||u_{n-1}|||_k^2) +C_2.\end{equation}

Arguing as in the proof of Lemma 4.3 and applying Lemmas 2.9 and 2.10 to inequality (4.26), we obtain
$|||u_n|||_k^2\leq C(k)$ for   $t\in [0,t_1]$. This proves the desired estimate.

\end{proof}

\subsection{Estimate for $\int\limits_{-\infty}^{+\infty} x^{2m} u^2_n dx$.}

\begin{lem} \label{estimate x_mu}
    Let $t_1$ be same as in Lemma \ref{norma2}. Then 
    $$\int\limits_{-\infty}^{+\infty} x^{2m} u^2_n dx\leq M_1(m),$$
    for all $0\leq t\leq t_1$, where $m\geq 1$ and $M_1(m)$ is a positive constant.
\end{lem}
\begin{proof}

We do this in similar way as in the case of previous estimates. First, we consider the case $m\geq 2$. So, we multiply each term of the equation by $2 x^{2m} u_n $ and integrate with respect to $x$ in $\mathbb R$.
\begin{equation}\label{ESTx^2m}
2\int\limits_{-\infty}^{+\infty} x^{2m} u_n \partial^\beta _{0t} u_n dx+2\int\limits_{-\infty}^{+\infty} x^{2m} u_n Lu_{n}dx
=2\gamma\int\limits_{-\infty}^{+\infty} x^{2m} u_n u^p_{n-1} u_{n-1,x} dx+2\int\limits_{-\infty}^{+\infty} x^{2m} u_nf_n(t) g(x,t)dx.
\end{equation}

First, we estimate the first integral on the right-hand side of the equality \eqref{ESTx^2m}.

According to the estimates obtained above, we have 
\begin{equation}\label{modul ux}
\sup_{x\in\mathbb R}\left|\frac{\partial u_{n-1}}{\partial x}\right|\leq\sqrt2\left(\int\limits_{-\infty}^{+\infty}\left(\frac{\partial u_{n-1}}{\partial x}\right)^2 dx \int\limits_{-\infty}^{+\infty}\left(\frac{\partial^2 u_{n-1}}{\partial x^2}\right)^2 dx\right)^\frac{1}{4}<C,
\end{equation}

and then

$$\bigg|2\gamma\int\limits_{-\infty}^{+\infty} u_n u^p_{n-1} \frac{\partial u_{n-1}}{\partial x} x^{2m} dx\bigg|\leq \sup_{x\in \mathbb R}\left|\frac{\partial u_{n-1}}{\partial x}\right|\sup_{x\in \mathbb R}(|u_{n-1}|)^{p-1} \int\limits_{-\infty}^{+\infty}|u_n||u_{n-1}|x^{2m} dx$$
\begin{equation}
\leq C\bigg(\int\limits_{-\infty}^{+\infty}x^{2m}u^2_ndx+\int\limits_{-\infty}^{+\infty}x^{2m}u^2_{n-1}dx\bigg).
\end{equation}

  For the second term on the right-hand side of equality \eqref{ESTx^2m}, taking into account estimates\eqref{estimate f_n} and \eqref{max f_n}, we have

\begin{equation}
\bigg|2\int\limits_{-\infty}^{+\infty} f_n(t)g(x,t) x^{2m} u_n dx \bigg|\leq C\int\limits_{-\infty}^{+\infty} x^{4m}g^2(x,t)  dx+\int\limits_{-\infty}^{+\infty} u_n^{2} dx\leq C_1,\end{equation}

Now we estimate integrals on the left-hand side of the equality \eqref{ESTx^2m}. According to the Lemma \ref{LemmAli} we have

\begin{equation}
2\int\limits_{-\infty}^{+\infty} x^{2m} u_n \partial^\beta_{0t} u_n dx\geq \partial^\beta_{0t} \int\limits_{-\infty}^{+\infty} x^{2m} u^2_n dx,\end{equation}

We simplify the other terms on the left-hand side by successively applying integration by parts.
$$2\int\limits_{-\infty}^{+\infty} x^{2m} u_n u_{nxxxx}dx=2 m(2m-1)(2m-2)(2m-3)\int\limits_{-\infty}^{+\infty} x^{2m-4} u^2_n dx $$

\begin{equation}
-8 m (2m-1)\int\limits_{-\infty}^{+\infty} x^{2m-2}u^2_{nx}dx+2\int\limits_{-\infty}^{+\infty} x^{2m} u^2_{nxx}dx.\end{equation}

Similarly, 
\begin{equation}
\int\limits_{-\infty}^{+\infty} x^{2m} u_n u_{nxxx} dx=M_1\int\limits_{-\infty}^{+\infty} x^{2m-3} u^2_n dx
+M_2\int\limits_{-\infty}^{+\infty} x^{2m-1} u^2_{nx} dx,\end{equation}

\begin{equation}
\int\limits_{-\infty}^{+\infty} x^{2m} u_n u_{nxx} dx=-\int\limits_{-\infty}^{+\infty} x^{2m} u^2_{nx} dx
+M_3\int\limits_{-\infty}^{+\infty} x^{2m-2} u^2_{n} dx,\end{equation}

\begin{equation}
2\int\limits_{-\infty}^{+\infty} x^{2m} u_n u_{nx} dx=-2m\int\limits_{-\infty}^{+\infty} x^{2m-1} u^2_n dx,\end{equation}

Now we need some inequalities for further estimations. For $v\in S(\mathbb R)$ and $k\in \mathbb N$ we have 

\begin{equation}
\left|\int\limits_{-\infty}^{+\infty} x^{k} v^2 dx\right|\leq \int\limits_{|x|\geq 1} |x|^{k} v^2 dx+\int\limits_{-1}^{1} v^2 dx\leq \int\limits_{-\infty}^{+\infty} |x|^{k+1}v^2 dx+ \int\limits_{-\infty}^{+\infty}v^2 dx, \end{equation}
and
$$
\left|\int\limits_{-\infty}^{+\infty} x^{k} v^2_x dx\right|= \left|-\int\limits_{-\infty}^{+\infty} x^{k} v v_{xx} dx+\frac{k(k-1)}{2}\int\limits_{-\infty}^{+\infty}x^{k-2}v^2 dx\right|\leq$$

$$\leq \sigma\int\limits_{-\infty}^{+\infty} |x|^{k} v_{xx}^2 dx +C(k,\sigma)\int\limits_{-\infty}^{+\infty}(|x|^{k}+1)v^2 dx.$$

So, taking in mind the above estimations and choosing appropriate value of $\sigma$, from \eqref{ESTx^2m} we obtain

\begin{equation}\label{BAHO}
\partial^\beta_{0t}\int\limits_{-\infty}^{+\infty} x^{2m} u^2_n
dx\leq M_5(m)\int\limits_{-\infty}^{+\infty} x^{2m} u^2_n dx+M_6
\int\limits_{-\infty}^{+\infty} x^{2m} u^{2}_{n-1} dx+M_7(m).\end{equation}
Here $M_6$ does not depend on $m$.

Applying the fractional Gronwall inequality  from  Lemma \ref{GRONWALL} to \eqref{BAHO}, we obtain 
$$\int\limits_{-\infty}^{+\infty} x^{2m} u^2_n dx\leq E_\beta(M_5T^\beta) \int\limits_{-\infty}^{+\infty} x^{2m}\varphi^2(x)dx $$

$$+
\frac{M_7T^\beta}{\beta} E_{\beta,\beta}(M_5T^\beta) + M_6 E_{\beta,\beta}(M_5T^\beta) I_{0t}^\beta \int\limits_{-\infty}^{+\infty} x^{2m} u^2_{n-1} dx=M_8+M_9I_{0t}^\beta\int\limits_{-\infty}^{+\infty} x^{2m} u^2_{n-1} dx.$$

 Applying Lemma \ref{yangi lemma1} we get
$$\int\limits_{-\infty}^{+\infty} x^{2m} u^2_n dx\leq M_8 E_\beta(M_9T^\beta)+M_9^nI_{0t}^{n\beta}\int\limits_{-\infty}^{+\infty} x^{2m} \varphi^2(x) dx$$

$$=M_8 E_\beta(M_9T^\beta)+\frac{M_9^n T^{n\beta}}{\Gamma(n\beta+1)}\int\limits_{-\infty}^{+\infty} x^{2m} \varphi^2(x) dx \leq const,$$
as $$\lim_{n\to +\infty} \frac{M_9^n T^{n\beta}}{\Gamma(n\beta+1)}=0.$$

Finally, the case $m=1$ follows from the case $m=2$ and Lemma 4.3. Indeed,
\[
\begin{aligned}
\int_{\mathbb R}x^2u_n^2\,dx
&=
\int_{|x|\leq 1}x^2u_n^2\,dx
+\int_{|x|>1}x^2u_n^2\,dx \\
&\leq
\int_{\mathbb R}u_n^2\,dx
+\int_{\mathbb R}x^4u_n^2\,dx
\leq C.
\end{aligned}
\]
Thus, the estimate also holds for $m=1$.
\end{proof}

A-priori estimates with respect to other semi-norms \eqref{yarim norma} can be obtained using equation \eqref{chiziqli tenglama} and inequality from Lemma \ref{Lemm_yakupov}. 

Now, to complete the proof of the Theorem \ref{th4.1}, we should show that $\{u_n\}$ and $\{f_n\}$ are Cauchy sequences for each of the above  semi-norms \eqref{yarim norma}.

\subsection{Convergence of $\{u_n\}$ and $\{f_n\}$}

Denote the difference $\omega_n=u_n-u_{n-1}$, $n\geq 1.$

Then for $n\geq 2$ we have the following equation for $\omega_n$.\begin{equation}\label{difference}
    \partial _{0t}^{\beta }{{\omega }_{n}}+L\omega_{n}
=\gamma u^p_{n-1}(u_{n-1,x}-u_{n-2,x})+\gamma u_{n-2,x}(u^p_{n-1}-u^p_{n-2})+(f_n(t)-f_{n-1}(t))g(x,t),\end{equation}
\begin{equation}
\omega_n(x,0)=0.
\end{equation}
Now we multiply the  equation \eqref{difference} by $\omega_n$ and integrate over $\mathbb{R}$.
$$\int\limits_{-\infty}^{+\infty}{\partial _{0t}^{\beta }{{\omega }_{n}}{{\omega }_{n}}dx
+\int\limits_{-\infty}^{+\infty}L\omega_{n}}{{\omega }_{n}}dx
=\gamma\int\limits_{-\infty}^{+\infty}u^p_{n-1}\omega_{n-1,x}\omega_n dx+\gamma\int\limits_{-\infty}^{+\infty}u_{n-2,x}(u^p_{n-1}-u^p_{n-2})\omega_ndx$$
\begin{equation} \label{convergence u_n f_n}
+\int\limits_{-\infty}^{+\infty}(f_n-f_{n-1})g\omega_ndx =I_1+I_2+I_3.\end{equation}
We estimate the right-hand side of the equality \eqref{convergence u_n f_n}. To do this, we use integration by parts and apply the Cauchy–Bunyakovsky–Schwarz inequality.

In the first term, using the estimate \eqref{modul ux}, we get
 
$$I_1=\gamma\int\limits_{-\infty}^{+\infty} u^p_{n-1}\omega_n\omega_{n-1,x}dx=-\gamma\int\limits_{-\infty}^{+\infty}( u^p_{n-1}\omega_n)_x\omega_{n-1}dx=$$
$$=-\gamma p\int\limits_{-\infty}^{+\infty}u^{p-1}_{n-1}u_{n-1,x}\omega_n\omega_{n-1}dx-\gamma\int\limits_{-\infty}^{+\infty}u^p_{n-1}\omega_{nx}\omega_{n-1}dx\equiv I_{1,1}+I_{1,2}.$$     
$$|I_{1,1}|\leq |\gamma||p|\sup_{x\in \mathbb R}|u^{p-1}_{n-1}|\sup_{x\in \mathbb  R}|u_{n-1,x}|\int\limits_{-\infty}^{+\infty}|\omega_n||\omega_{n-1}|dx\leq|\gamma||p|CM\bigg( \frac{1}{2}\|\omega_n\|^2+\frac{1}{2}\|\omega_{n-1}\|^2\bigg),$$

$$|I_{1,2}|\leq |\gamma|\sup_{x\in \mathbb R}|u^p_{n-1}|\int\limits_{-\infty}^{+\infty}|\omega_{nx}||\omega_{n-1}|dx\leq C\bigg( \|\omega_{nx}\|^2+\|\omega_{n-1}\|^2\bigg) $$

$$\leq \varepsilon\|\omega_{nxx}\|^2+C(\varepsilon)\bigg(\|\omega_n\|^2+\|\omega_{n-1}\|^2\bigg).$$

Now, we estimate integral
$$I_2=\gamma\int\limits_{-\infty}^{+\infty}u_{n-2,x}(u^p_{n-1}-u^p_{n-2})\omega_n dx=\gamma\int\limits_{-\infty}^{+\infty}\omega_nu_{n-2,x}(u_{n-1}-u_{n-2})\sum^{p-1}_{k=0}u^{p-1-k}_{n-1}u^k_{n-2}dx.$$

$$|I_2|\leq |\gamma|\sup_{x\in \mathbb R }|u_{n-2,x}|\sup_{x\in \mathbb R }\bigg|\sum^{p-1}_{k=0}u^{p-1-k}_{n-1}u^k_{n-2}\bigg|\|\omega_{n-1}\|\|\omega_n\|.$$

Taking into account the inequality $ |v|\leq\sqrt{2} \|v\|\|v_x\|$ for any $v\in S(\mathbb R)$ and using Lemmas \ref{norma2}, \ref{high derivatives} we get 
\begin{equation}\label{I^3}
|I_2|\leq C\bigg(\|\omega_n\|^2+\|\omega_{n-1}\|^2\bigg).\end{equation}

Now we estimate $|I_3|$.
$$f_n(t)-f_{n-1}(t)=(Bf_n)(t)-(Bf_{n-1})(t)=$$
$$=\frac{1}{g^*(t)}\bigg[a^2\int\limits_{-\infty}^{+\infty}\eta_{xxxx}(x)\omega_{n}dx-b\int\limits_{-\infty}^{+\infty}\eta_{xxx}(x)\omega_{n}dx$$
$$+c\int\limits_{-\infty}^{+\infty}\eta_{xx}(x)\omega_{n}dx-d\int\limits_{-\infty}^{+\infty}\eta_x(x)\omega_{n}dx+q\int\limits_{-\infty}^{+\infty}\eta(x)\omega_ndx+\frac{\gamma}{p+1}\int\limits_{-\infty}^{+\infty}\eta_x(x)\omega_{n-1}\sum_{j=0}^{p}u^{p-j}_{n-1}u^j_{n-2}dx\bigg]$$

$$=\frac{1}{g^*(t)}(I_{3,1}+I_{3,2}+I_{3,3}+I_{3,4}+I_{3,5}+I_{3,6})$$
$$|I_{3,1}|\leq a^2\|\eta_{xxxx}\|_{L^2}\|\omega_{n}\|_{L^2}, \quad|I_{3,2}|\leq |b|\|\eta_{xxx}\|_{L^2}\|\omega_{n}\|_{L^2},$$
$$|I_{3,3}|\leq|c|\|\eta_{xx}\|_{L^2}\|\omega_{n}\|_{L^2},\quad |I_{3,4}|\leq |d|\|\eta_x\|_{L^2}\|\omega_{n}\|_{L^2},\quad |I_{3,5}|\leq |q|\|\eta\|_{L^2}\|\omega_n\|_{L^2}.$$

As
$$\bigg|\sum^p_{j=0}u^{p-j}_{n-1}u^j_{n-2}\bigg|\leq const.$$
we have 
$$|I_{3,6}|\leq |\gamma|M^p\int\limits_{-\infty}^{+\infty}|\eta_x||\omega_{n-1}|dx\leq C_2\|\omega_{n-1}\|_{L^2}.$$
Thus, by combining all terms, we obtain the following inequality
\begin{equation}\label{dsd}
|f_n-f_{n-1}|\leq C_4\|\omega_n\|+C_5\|\omega_{n-1}\|.\end{equation}
Substituting the inequality \eqref{dsd} into the estimate \eqref{convergence u_n f_n}, we obtain the following
\begin{equation}
|I_3|\leq \bigg|(f_n-f_{n-1})\int\limits_{-\infty}^{+\infty}g\omega_ndx\bigg|\leq \|g\|_{L^2}C\bigg(\|\omega_n\|+\|\omega_{n-1}\|\bigg)\|\omega_n\|\leq C_6\|\omega_n\|^2+C_7\|\omega_{n-1}\|^2.\end{equation}

Taking the above estimates for $I_1$, $I_2$, and  $I_3$ into account, we have 
$$|I|\leq C(\|\omega_n\|^2+\|\omega_{n-1}\|^2).$$

Choosing appropriate $\varepsilon$, in the same manner as in above obtained estimates, we get 

$$\partial _{0t}^{\beta }{{\left\| {{\omega }_{n}} \right\|}^{2}}\le \widetilde{C}\left[ {{\left\| {{\omega }_{n}} \right\|}^{2}}+{{\left\| {{\omega }_{n-1}} \right\|}^{2}} \right].$$


From the last inequality, using the Gronwall's inequality from Lemma \ref{GRONWALL}, we get

\begin{equation}
{{\left\| {{\omega }_{n}}(\cdot,t) \right\|}^{2}} \leq {{\left\| {{\omega }_{n}}(\cdot,0) \right\|}^{2}}E_\beta(\widetilde{C}t^\beta)+\widetilde{C} \cdot\Gamma(\beta)E_{\beta,\beta}(\widetilde{C}T^\beta)I^\beta _{0t}{{\left\| {{\omega }_{n-1}}(\cdot,t) \right\|}^{2}} = KI^\beta _{0t}{{\left\| {{\omega }_{n-1}}(\cdot,t) \right\|}^{2}}.
\end{equation}

Here we used the fact that $\omega_n(x,0)=u_n(x,0)-u_{n-1} (x,0)=0$.

Therefore, based on \ref{yangi lemma1}, the following inequality is obtained:

\begin{equation}
{{\left\| {{\omega }_{n}}(\cdot,t) \right\|}^{2}} \leq K^{n-1}I_{0t}^{(n-1)\beta}{{\left\| {{\omega }_{1}}(\cdot,t) \right\|}^{2}}.\end{equation}

From the above estimates , since
$\|\omega _1 (\cdot,t)\|=\|u _1 (\cdot,t)-u_0(\cdot,t)\|\leq const$
we have
\begin{equation}
{{\left\| {{\omega }_{n}}(\cdot,t) \right\|}^{2}} \leq K^{n-1}C\frac{T^{(n-1)\beta}}{\Gamma\left((n-1)\beta+1\right)}.\end{equation}

Now we prove that the sequence $\{u_n\}$ is a Cauchy (fundamental) sequence with respect to norm $\|\cdot\|$.
We have

$$\|u_{n+p} (\cdot,t)-u_n (\cdot,t)\|=\left\|\sum_{k={n+1}}^{n+p}\omega_k (\cdot,t)\right\|\leq\sum_{k={n+1}}^{n+p}\|\omega_k (\cdot,t)\|\leq C\sum_{k=n+1}^{n+p}\frac{(KT^\beta)^\frac{{k-1}}{2}}{\sqrt {\Gamma((k-1)\beta+1)}}.$$

It is well known that the series 
$\sum_{k=1}^{\infty}\frac{(KT^\beta)^\frac{{k-1}}{2}}{\sqrt {\Gamma((k-1)\beta+1)}}$
is a convergent series. Therefore, by the Cauchy criterion,

$$\forall\varepsilon>0,\exists N: \forall n\geq N, \ \ \sum_{k=n+1}^{n+p}\|\omega_k(\cdot,t)\|<\varepsilon.$$
Consequently,
$$\max_t\left\|u_{n+p} (\cdot,t)-u_n (\cdot,t)\right\|<\varepsilon.$$
This means that the sequence $u_n (x,t)$ is a Cauchy (fundamental) sequence in the norm of $C([0,t_1];L^2(\mathbb R))$. 

Now we prove that $f_n(t)$ is a Cauchy sequence. Based on \eqref{convergence u_n f_n},
$$|f_{n+p}-f_n|\leq \sum_{k=n+1}^{n+p}|f_k-f_{k-1}|\leq C_4\sum_{k=n+1}^{n+p}\|\omega _k\|+C_5\sum_{k=n+1}^{n+p}\|\omega_{k-1}\|=C_4\sum_{k=n+1}^{n+p}\|\omega_k\|+C_5\sum_{k=n}^{n+p-1}\|\omega_k\|$$
$$=C_5\|\omega_n\|+(C_4+C_5)\sum_{k=n+1}^{n+p-1}\|\omega_k\|+C_4\|\omega_{n+p}\|\leq M(\|\omega_n\|+2\sum_{k=n+1}^{n+p-1}\|\omega_k\|+\|\omega_{n+p}\|)$$
$$\leq2M\sum_{k=n}^{n+p}\|\omega_k\|\leq 2M \sqrt{C}\sum_{k=n}^{n+p}\frac{(KT^\beta)^{\frac{k-1}{2}}}{\sqrt{\Gamma((k-1)\beta+1)}}\leq A\sum_{k=n}^{n+p}\frac{(KT^\beta)^{\frac{k-1}{2}}}{\sqrt{\Gamma((k-1)\beta+1)}},$$

$$\lim_{n\rightarrow{+\infty}}A\sum_{k=n}^{n+p}\frac{(KT^\beta)^{\frac{k-1}{2}}}{\sqrt{\Gamma((k-1)\beta+1)}}=0.$$
 Therefore,
 $$\forall\varepsilon>0,\exists N: \forall n\geq N, \ \ |f_{n+p}-f_n|<\varepsilon,$$ hence $\{f_n\}$ is a Cauchy sequence.

For the semi-norms involving higher-order derivatives, we use the Gagliardo-Nirenberg inequality as follows. 
By the Gagliardo-Nirenberg inequality Lemma \ref{Gagliardo ineq} with 
$d=1, j=k, m=k+1, q=r=s=2$ we have  $\theta=\frac{k}{k+1}$ and therefore 
$$\|u_{n+p}-u_n\|_k\leq C\|u_{n+p}-u_n\|^{\frac{k}{k+1}}_{k+1}\|u_{n+p}-u_n\|^{\frac{1}{k+1}}\leq \tilde C\|u_{n+p}-u_n\|^{\frac{1}{k+1}}.$$
Since, the sequence is Cauchy in the semi-norms $|||\cdot |||_k$.

Now we estimate $||{{x}^{m}}(u_{n+p}-u_n)||^2=\int\limits_{-\infty }^{+\infty }x^{2m}{(u_{n+p}-u_n)^{2}}dx$. 
Using the Holder's inequality we have
$$\int\limits_{-\infty }^{+\infty }x^{2m}{(u_{n+p}-u_n)^{2}}dx \leq \left(\int\limits_{-\infty }^{+\infty }{(u_{n+p}-u_n)^{2}}dx\int\limits_{-\infty }^{+\infty }x^{4m}{(u_{n+p}-u_n)^{2}}dx\right)^{1/2}\leq C \|u_{n+p}-u_n\|.$$
From the last inequality follows that $\{u_n\}$ is Cauchy in weighted semi-norms.





Based on the above estimates, we conclude that the sequence
$\{u_n\}$ converges in $SC(t_1)$ to some function $u$ and that
$f_n\to f$ in $C[0,t_1]$. From equation (4.4), we have
\[
\partial_{0t}^{\beta}u_n
=-Lu_n+f_n(t)g+\gamma u_{n-1}^{p}u_{n-1,x}.
\]
Since $u_n\to u$ in $SC(t_1)$ and $f_n\to f$ in $C[0,t_1]$, the
right-hand side converges in $SC(t_1)$ to
\[
-Lu+f(t)g+\gamma u^{p}u_x.
\]
By the closed-graph consequence of the fundamental theorem of
fractional calculus~[34, Theorem~3.7], it follows that
\[
\partial_{0t}^{\beta}u
=-Lu+f(t)g+\gamma u^{p}u_x\in SC(t_1).
\]
Therefore, $u\in SC_{\beta}(t_1)$, and passing to the limit in
(4.4) shows that $(u,f)$ is a solution of the inverse problem.
This completes the proof of the theorem.

\section{CONTINUOUS DEPENDENCE ON THE GIVEN DATA}
\begin{thm} \label{Stability}
Let $\varphi,$ $\tilde \varphi\in S(\mathbb R)$ and $E,$ $\tilde E$, $\partial^\beta_{0t}E$, $\partial^\beta_{0t}\tilde E\in C[0,T].$ Suppose that $(u,f)$ and $(\tilde u,\tilde f)$ are the corresponding solutions to the inverse problem \eqref{41}--\eqref{41-addC} associated with initial and additional data $\varphi,E,$ and $\tilde \varphi,\tilde E$, respectively, on common time interval $[0,\tau]$, $0<\tau\leq T$. Then  for every pair of nonnegative integers $m$ and $k$, the following stability estimates hold.

$$\max_{0\leq t\leq \tau}||u-\tilde u||^2_{H^k(\mathbb R)}+\max_{0\leq t\leq \tau}\|x^m (u-\tilde u)\|^2\leq M_{m,k}(||\varphi-\tilde\varphi||^2_{H^k(\mathbb R)}+||x^m(\varphi-\tilde\varphi)||^2+\max_{0\leq t\leq \tau}|\partial^\beta_{0t}(E-\tilde E)|^2), $$

$$\max_{0\leq t\leq \tau}|f-\tilde f|\leq M_*\bigg(||\varphi-\tilde \varphi||+\max_{0\leq t\leq \tau}|\partial^\beta_{0t}(E(t)-\tilde E(t))|\bigg ),$$
where $M_{m,k}$ and $M_*$ are positive constants. 
\end{thm}

\begin{proof}
We denote
$v=u-\tilde u$, and we get the following equation for $v$
\begin{equation}\label{stability equation}
\partial _{0t}^{\beta }{v}+Lv=\gamma (u^pu_x-\tilde u^p\tilde u_x)+(f-\tilde f)g. \quad x\in \mathbb R, \ 0<t\leq T.\end{equation}
with initial condition $$v(x,0)=\varphi(x)-\tilde\varphi(x), \quad x\in \mathbb R.$$
and additional condition
$$\int\limits_{-\infty}^{+\infty}\eta(x)v(x,t)dx=E(t)-\tilde E(t).$$
We multiply both sides of this equation by $2v$ and integrate in $\mathbb R$ with respect to $x$.
$$2\int\limits_{-\infty}^{+\infty}v\partial^\beta_{0t}vdx+2a^2\int\limits_{-\infty}^{+\infty}vv_{xxxx}dx+2b\int\limits_{-\infty}^{+\infty}vv_{xxx}dx+2c\int\limits_{-\infty}^{+\infty}vv_{xx}dx+2d\int\limits_{-\infty}^{+\infty}vv_xdx$$
\begin{equation}
+2q\int\limits_{-\infty}^{+\infty}v^2dx=2\int\limits_{-\infty}^{+\infty}v(f-\tilde f)gdx+2\gamma \int\limits_{-\infty}^{+\infty}(u^pu_x-\tilde u^p\tilde u_x)vdx.\end{equation}
The terms in the left hand side can be estimated in same manner as in above paragraphs. We now estimate the right-hand side terms. 

Repeating similar estimation as in the proof of inequality \eqref{I^3} we get 

$$ \bigg| \int\limits_{-\infty}^{+\infty}(u^pu_x-\tilde u^p\tilde u_x)vdx\bigg|\leq  C\|v\|^2$$
and using \eqref{B_n}
\begin{equation} \label{f-f}
    |f-\tilde f|\leq C\left(\|v\|+|\partial_{0t}^\beta (E-\tilde E)|\right).
\end{equation}
So, applying the Cauchy-Schwarz inequality, we obtain 
$$\bigg|2(f-\tilde f)\int\limits_{-\infty}^{+\infty}gvdx\bigg|\leq \|v\|^2+\|g\|^2(f-\tilde f)^2\leq (1+C^2\|g\|^2)\|v\|^2+\frac{\|g\|^2}{(\min_t|g^*(t)|)^2}[\partial^\beta_{0t}(E(t)-\tilde E(t))]^2.$$

Summarizing, we get
$$\partial^\beta_{0t}\|v\|^2\leq C_2\|v\|^2+C_3[\partial^\beta_{0t}(E(t)-\tilde E(t))]^2.$$

Now, applying Gronwall's inequality \eqref{GRONWALL} to the resulting inequality, we finally obtain the first stability estimate:

$$\sup_t\|u(\cdot,t)-\tilde u(\cdot,t)\|^2\leq M\bigg(\|\varphi-\tilde \varphi\|^2+\sup_t[\partial^\beta_{0t}(E(t)-\tilde E(t))]^2\bigg ).$$

This inequality establishes the stability of the solution in $C([0,T] ; L_2(\mathbb R))$ norm, i.e. a small perturbation in the initial data and in the parameter leads to a correspondingly small deviation in the solution. From \eqref{f-f} we get
$$ \max_t |f-\tilde f|\leq C\|v\|\leq M_*\bigg(\|\varphi-\tilde \varphi\|+\sup_t\partial^\beta_{0t}(E(t)-\tilde E(t))|\bigg ).$$

Take  $k$-th  -derivative on each sides of \ref{stability equation}, multiply by $\frac{{{\partial }^{k}}v}{\partial {{x}^{k}}}$ and integrate over $\mathbb R$. Estimating the linear terms as above, we obtain

$$
\begin{aligned}
\frac{1}{2}\partial_{0t}^{\beta}|||v|||_k^2
+a^2|||v|||_{k+2}^2
&\leq
C\left(|||v|||_k^2+|||v|||_{k+1}^2\right)
+C|f-\widetilde f|^2
\\
&\quad
+|\gamma|
\left|
\int_{\mathbb R}
\partial_x^k
\left(u^pu_x-\widetilde u^{\,p}\widetilde u_x\right)
\partial_x^kv\,dx
\right|.
\end{aligned}
$$

To estimate the nonlinear term, we write
\[
u^pu_x-\widetilde u^{\,p}\widetilde u_x
=\frac{1}{p+1}\partial_x
\left(u^{p+1}-\widetilde u^{\,p+1}\right)
=\frac{1}{p+1}\partial_x
\left(vP_p(u,\widetilde u)\right),
\]
where
\[
P_p(u,\widetilde u)
=\sum_{j=0}^{p}u^{p-j}\widetilde u^{\,j}.
\]
By the Leibniz formula and the uniform estimates for $u$ and
$\widetilde u$ obtained in Lemma~4.4, we have
\[
\begin{aligned}
&\left|
\int_{\mathbb R}
\partial_x^k
\left(u^pu_x-\widetilde u^{\,p}\widetilde u_x\right)
\partial_x^kv\,dx
\right|\leq
C_k\left(|||v|||_{k+1}+\|v\|_{H^k(\mathbb R)}\right)|||v|||_k\leq
\\
&\qquad\leq
\varepsilon|||v|||_{k+2}^2
+C_{k,\varepsilon}\|v\|_{H^k(\mathbb R)}^2.
\end{aligned}
\]
Here we used
\[
|||v|||_{k+1}^2
\leq
|||v|||_k|||v|||_{k+2}
\leq
\varepsilon|||v|||_{k+2}^2
+\frac{1}{4\varepsilon}|||v|||_k^2.
\]
Moreover, it follows from (5.3) that
\[
|f-\widetilde f|^2
\leq
C\left(
\|v\|^2
+\left|\partial_{0t}^{\beta}(E-\widetilde E)\right|^2
\right).
\]
Therefore, choosing $\varepsilon>0$ sufficiently small and summing
the resulting estimates for the derivatives of orders $0,\ldots,k$,
we obtain
\[
\partial_{0t}^{\beta}\|v\|_{H^k(\mathbb R)}^2
\leq
C_k\left(
\|v\|_{H^k(\mathbb R)}^2
+\left|\partial_{0t}^{\beta}(E-\widetilde E)\right|^2
\right).
\]

Now we estimate $||{{x}^{m}}v||^2=\int\limits_{-\infty }^{+\infty }{{{x}^{2m}}{{v}^{2}}dx}$.  For that we mostly use the same method as in the proof of Lemma \ref{estimate x_mu}. We have
	\[2\int\limits_{-\infty }^{+\infty }{{{x}^{2m}}}v\partial _{0t}^{\beta }vdx+2\int\limits_{-\infty }^{+\infty }{{{x}^{2m}}}vLvdx=2\gamma \int\limits_{-\infty }^{+\infty }{{{x}^{2m}}}\left( {{u}^{p}}{{u}_{x}}-{{{\tilde{u}}}^{p}}{{{\tilde{u}}}_{x}} \right)vdx+2\int\limits_{-\infty }^{+\infty }{{{x}^{2m}}}g(f-\tilde{f})vdx.\]

    First, we estimate second term in the right-hand side of this inequality. 
\[\left| \int\limits_{-\infty }^{+\infty }{{{x}^{2m}}}g(f-\tilde{f})vdx \right|\le |f-\tilde{f}|\,\,||{{x}^{2m}}g||\,\,||v||\le C\left( ||\varphi -\tilde{\varphi }|{{|}^{2}}+\underset{t}{\mathop{\sup }}\,|\partial _{0t}^{\beta }(E-\tilde{E}){{|}^{2}} \right).\]

Now we estimate nonlinear term.
	\[\left| \int\limits_{-\infty }^{+\infty }{{{x}^{2m}}}\left( {{u}^{p}}{{u}_{x}}-{{{\tilde{u}}}^{p}}{{{\tilde{u}}}_{x}} \right)vdx \right|=\left| \int\limits_{-\infty }^{+\infty }{{{x}^{2m}}}\left( {{u}^{p}}-{{{\tilde{u}}}^{p}} \right){{u}_{x}}vdx \right|+\left| \int\limits_{-\infty }^{+\infty }{{{x}^{2m}}}{{{\tilde{u}}}^{p}}{{v}_{x}}vdx \right|\]
	\[\le \left( \underset{x}{\mathop{\sup }}\,\left| {{x}^{2m}}{{u}_{x}}\sum\limits_{k=0}^{p-1}{{{u}^{k}}{{u}^{p-1-k}}} \right|+\frac{1}{2}\underset{x}{\mathop{\sup }}\,\left| {{\left( {{x}^{2m}}{{{\tilde{u}}}^{p}} \right)}_{x}} \right| \right){{\left\| v \right\|}^{2}}.\]
The expression in the brackets is bounded due to inequality \eqref{SN-SN} and Lemmas \ref{norma2}---\ref{estimate x_mu}.

So, we get
$$\left| \int\limits_{-\infty }^{+\infty }{{{x}^{2m}}}\left( {{u}^{p}}{{u}_{x}}-{{{\tilde{u}}}^{p}}{{{\tilde{u}}}_{x}} \right)vdx \right|\leq C\left( ||\varphi -\tilde{\varphi }|{{|}^{2}}+\max_{0\leq t\leq\tau}|\partial _{0t}^{\beta }(E-\tilde{E}){{|}^{2}} \right).$$

Other terms can be estimated as in the proof of Lemma \ref{estimate x_mu}. So, we have 

$$\partial _{0t}^{\beta }\|x^mv\|^2\leq C \|x^mv\|^2+  C\left( ||\varphi -\tilde{\varphi }|{{|}^{2}}+\max_{0\leq t\leq\tau}|\partial _{0t}^{\beta }(E-\tilde{E}){{|}^{2}} \right).$$

 Using Gronwall-Bellman's inequality, we obtain
$$\|x^mv\|^2\leq  C\left( ||x^m(\varphi -\tilde{\varphi })|{{|}^{2}}+||\varphi -\tilde{\varphi }|{{|}^{2}}+\max_{0\leq t\leq\tau}|\partial _{0t}^{\beta }(E-\tilde{E})|^{2} \right). $$

From the above obtained estimates we get the stability estimates in the Theorem. This completes the proof of the theorem.
\end{proof}

\textbf{Remark.} Uniqueness of the solution follows from the stability estimates.

\section{Finite-time blow-up of solutions}

Now, let us consider the case when the nondegeneracy quantity $g^{*}(t)$ is a
strictly positive constant. In this case one might expect the inverse
problem to be globally well posed. We show, however, that this is not so:
even under the strongest possible form of the nondegeneracy condition, the
solution may blow up in finite time. Consider the following example.

\medskip
\noindent\textbf{Example.} Let $T^{*}>0$. We consider an inverse source problem 
\begin{equation}\label{example}
\begin{cases}
\partial_{0t}^{\beta}u+a^{2}u_{xxxx}=f(t)g(x,t)+\gamma u^{p}u_{x}, &
(x,t)\in\mathbb{R}\times(0,T^{*}),\\[4pt]
u(x,0)=(T^*)^{-\alpha}\psi(x), & x\in\mathbb{R},\\[4pt]
\int\limits_{-\infty}^{+\infty}\eta(x)u(x,t)\,dx=0, & t\in[0,T^{*}),
\end{cases}\end{equation}
where $0<\beta<1$, $p\in\mathbb{N}$, $a>0$, $\gamma>0$, $\alpha=\frac{\beta}{p}\,$, 
$
 \psi(x)=e^{-x^{2}},\,\, \eta(x)=xe^{-x^{2}}.$ We define  
 $$g(x,t)=h(t)\psi(x)+a^{2}r(t)\psi^{(4)}(x)-\gamma\,\psi^{p}(x)\psi'(x),
$$
where $r(t)=(T^{*}-t)^{\beta}$ and
$$
h(t)=\frac{\alpha}{\Gamma(1-\beta)}
\int_{0}^{t/(T^{*}-t)} z^{-\beta}(1+z)^{-\alpha-1}\,dz.
$$
for $t<T^*$.
The function $h(t)$ is nonnegative, nondecreasing, and continuous on
$[0,T^{*})$, and satisfies
$$
0\le h(t)\le \,\frac{\Gamma(\alpha+\beta)}{\Gamma(\alpha)},
\qquad
\lim_{t\to T^{*-}} h(t)=\frac{\Gamma(\alpha+\beta)}{\Gamma(\alpha)}.
$$
Consequently, if we extend $g$ for $t\geq T^*$ by setting $r(t)=0$ and $h(t)=\Gamma(\alpha+\beta)/\Gamma(\alpha)$, we get  $g\in C\big([0,\infty);S(\mathbb{R})\big)$.

Then $$
g^{*}(t)=\int\limits_{-\infty}^{+\infty}\eta(x)g(x,t)\,dx
=-\gamma\int\limits_{-\infty}^{+\infty}\eta(x)\psi^{p}(x)\psi'(x)\,dx
=\frac{\gamma\sqrt{\pi}}{(p+2)^{3/2}}>0
\qquad\text{for all } t\ge0.
$$

One can verify that the pair of functions
$$
u(x,t)=(T^{*}-t)^{-\alpha}\psi(x), \qquad
f(t)=(T^{*}-t)^{-\alpha(p+1)}
$$
satisfies problem \eqref{example}, 
while at the same time
$$
\|u(\cdot,t)\|_{L^{2}(\mathbb{R})}+|f(t)|\longrightarrow+\infty
\qquad\text{as } t\to T^{*-},
$$
for all $0<\beta<1$ and $p\in\mathbb{N}$.

This blow-up solutions shows that under the conditions of the Theorem \ref{th4.1} the global solvability could not be proven.

\section{CONCLUSION}

In this paper, an inverse source problem for a generalized
time-fractional Kuramoto--Sivashinsky equation under an integral
overdetermination condition was investigated. For the linear problem,
existence, uniqueness, and stability were established by reducing the
problem to an operator equation and using the quasinilpotency of the
corresponding operator. For the nonlinear problem, local existence in
the Schwartz space was proved by the method of successive
approximations. A priori estimates and the fractional Gronwall
inequality were used to prove uniqueness and derive stability
estimates for the solution and the source factor with respect to the
initial data and the Caputo derivative of the integral observation.

An explicit solution exhibiting finite-time blow-up was also
constructed. In this example, both
$\|u(\cdot,t)\|_{L_2(\mathbb R)}$ and $|f(t)|$ tend to infinity,
although $g^*(t)$ remains a strictly positive constant. This shows
that the nondegeneracy condition does not prevent finite-time
blow-up. Thus, the assumptions ensuring local solvability do not,
in general, guarantee global solvability of the nonlinear inverse
problem.

\label{lastpage}


%
\end{document}